\documentclass[draft, 10pt]{article}
\usepackage{amssymb, amsthm, amsmath}
\numberwithin{equation}{section}

\newcommand{\x}{\ma{x}}

\newcommand{\ph}{\mathsf{PHH}}
\newcommand{\mal}{\boldsymbol{\alpha}}
\newcommand{\cc}{$\mathsf{CC}$}
\newcommand{\cl}{$\mathsf{CL}$}

\newcommand{\R}{\mathbb{R}}

\newcommand{\Z}{\mathbb{Z}}
\newcommand{\N}{\mathbb{N}}
\newcommand{\Q}{\mathbb{Q}}
\newcommand{\Grass}{\mathbb{G}}
\newcommand{\bfP}{\mathbb{P}}

\newcommand{\e}{\emph}
\newcommand{\rom}{\mathrm}
\newcommand{\ov}{\overline}
\newcommand{\ma}{\mathbf}

\newcommand{\mcal}{\mathcal}

\newcommand{\ben}{\begin{enumerate}}
\newcommand{\een}{\end{enumerate}}
\newcommand{\eit}{\begin{itemize}}
\newcommand{\beq}{\begin{equation}}
\newcommand{\eeq}{\end{equation}}

\newcommand{\ve}{\varepsilon}

\newcommand{\al}{\alpha}

\newcommand{\del}{\delta}

\newcommand{\lab}{\label}

\newtheorem{thm}{Theorem}
\newtheorem*{thm*}{Theorem}
\newtheorem{lem}{Lemma}
\newtheorem{pro}{Proposition}

\newtheorem*{cor*}{Corollary}
\newtheorem*{hyp}{Hypothesis}
\newtheorem{con}{Conjecture}

\theoremstyle{definition}
\newtheorem*{ack}{Acknowledgements}

\newcommand{\barray}{\begin{eqnarray*}}
\newcommand{\earray}{\end{eqnarray*}}
\newcommand{\hcf}{\rom{h.c.f.}}

\newcommand{\colt}[2]{\genfrac{}{}{0pt}{1}{#1}{#2}}

\newtheorem{rmk}{Remark}
\newtheorem{prop}{Proposition}
\newtheorem{lemma}{Lemma}
\newtheorem*{coro}{Corollary}

\newcommand{\PP}{\mathbb{P}}

\newcommand{\lt}{\left}
\newcommand{\rt}{\right}

\newcommand{\OO}{\mathcal O}

\def\PP{{\mathbb P}}

\begin{document}

\title{The density of  rational points on\\ non-singular hypersurfaces, II}

\author{T.D. Browning$^1$ and D.R. Heath-Brown$^2$\\
\small{$^1$\emph{School of Mathematics,  
Bristol University, Bristol BS8 1TW}}\\
\small{$^2$\emph{Mathematical Institute,
24--29 St. Giles',
Oxford OX1 3LB}}\\
\small{$^1$t.d.browning@bristol.ac.uk}, \small{$^2$rhb@maths.ox.ac.uk}}
\date{}
\maketitle

\begin{abstract}
For any integers $d,n \geq 2$, let $X \subset \mathbb{P}^{n}$ be a non-singular hypersurface 
of degree $d$ that is defined over $\mathbb{Q}$. The main result in this paper
is a proof that the number $N_X(B)$ of $\mathbb{Q}$-rational points on $X$ which have
height at most $B$ satisfies
$$
N_X(B)=O_{d,\varepsilon,n}(B^{n-1+\varepsilon}),
$$
for any $\varepsilon>0$. The implied constant in this estimate depends
at most upon $d, \varepsilon$ and $n$.
\end{abstract}

\section{Introduction}

For any $n \geq 2$, let $F \in \Z[X_0,\ldots,X_n]$ be a form of degree
$d\geq 2$ that defines a non-singular hypersurface $X \subset \bfP^{n}$.
In this paper we return to the theme of our recent investigation
\cite{smoothI} into the distribution of rational points on such hypersurfaces.  
For any rational point 
$x=[\x] \in \bfP^n(\Q)$ such that
$\x=(x_0,\ldots,x_n) \in \Z^{n+1}$ and $\hcf(x_0,\ldots,x_n)=1$, we
shall write 
\beq\lab{height}
H(x)=|\x|
\eeq
for its height, where 
$|\x|$ denotes the norm $\max_{0\leq i \leq n}|x_i|$.
With this notation in mind, our primary objective is to understand the
asymptotic behaviour of the quantity
$$
N_X(B)=\#\{x \in X\cap\bfP^n(\Q): ~H(x) \leq B\},
$$
as $B\rightarrow \infty$.  We have the following basic conjecture.

\begin{con}\lab{hb-ns}  Let $\ve>0$.  Then we have
\beq\lab{n-1}
N_X(B)=O_{d,\ve,n}(B^{n-1+\ve}).
\eeq
\end{con}

Throughout our work the implied constant in any estimate is absolute unless 
explicitly indicated otherwise.  
In the case of Conjecture \ref{hb-ns}, for example, 
the constant is permitted to depend only upon $d$, $\ve$ and $n$. 
In view of the fact that $\x$ and $-\x$ represent the same point in
projective space, it is clear that 
$$
N_X(B)=\frac{1}{2}\#\{\x \in \Z^{n+1}: F(\x)=0, ~\hcf(x_0,\ldots,x_n)=1, ~|\x| \leq B\},
$$
and so one may equally view Conjecture \ref{hb-ns} as a statement
about the frequency of integer solutions to certain homogeneous Diophantine equations.
In particular, when $F$ is a non-singular quadratic form in $4$
variables, with discriminant equal to
a square, we have $N_X(B)=c_X B^{2}\log B (1+o(1))$. In all other
cases we would suppose that the exponent $\ve$ is
superfluous in Conjecture \ref{hb-ns}.  In fact there is a 
conjecture of Batyrev and Manin \cite{bat} that predicts
one should have $N_X(B)=O_{X}(B^{n-1-\delta})$, for some $\delta>0$, provided that $d \geq
3$ and $n \geq 4$.

Conjecture \ref{hb-ns} is a special case of a conjecture due
to the second author \cite[Conjecture 2]{annal}, which predicts that the same
estimate should hold under the weaker assumption that the defining
form $F$ is absolutely irreducible.
Both of these conjectures have received a significant amount of
attention in recent years, to the extent that Conjecture \ref{hb-ns} is now
known for all values of $d \geq 2$ and $n \geq 3$, except possibly for
the eight cases in which $d=3$ or $4$ and $n=5,6,7$ or $8$.
This comprises the combined work of both the first and second authors
\cite{3fold,bhb,smoothI,hb-india,annal}, as summarised in 
\cite[Corollary 1]{smoothI}.
For the exceptional cases the best result available is the estimate
$N_X(B) \ll_{d,\ve,n}B^{n-1+\theta_{d}+\ve}$, for any $\ve>0$, with
$$
\theta_{d}=\left\{
\begin{array}{ll}
5/(3\sqrt{3})-3/4, & d=3,\\
1/12, & d= 4.
\end{array}
\right.
$$
This follows from joint work of the authors with Salberger \cite{bhbs}.

The aim of the present paper is to complete the proof of Conjecture
\ref{hb-ns}, by offering a satisfactory treatment of the eight
remaining cases.  To this end we define the set
\beq\lab{E}
\mcal{E}=\{(d,n): \mbox{$d=3$ or $4$, ~$n=5,6,7$ or $8$}\}.
\eeq
The following is our primary result.

\begin{thm}\lab{6}
Let $\ve>0$ and $(d,n) \in \mcal{E}$, and suppose that $F \in \Z[X_0,\ldots,X_n]$
is a non-singular form of degree $d$.
Then we have
$$
N_X(B)=O_{d,\ve,n}(B^{n-1+\ve}).
$$
\end{thm}

\begin{cor*}
Conjecture \ref{hb-ns} holds in every case.
\end{cor*}

The authors have recently learnt of work due
to Salberger \cite{salberger}, which establishes Conjecture
\ref{hb-ns} in the case $d=4$.
In fact Salberger obtains the estimate (\ref{n-1}) 
for any geometrically integral hypersurface $X\subset \bfP^n$ of
degree $d\geq 4$, such that $X$ contains at most finitely many linear subspaces of dimension $n-2$.

Returning to the setting of non-singular hypersurfaces, 
the corollary to Theorem \ref{6} is in some sense most significant
when $d \leq n$, since it is precisely in this setting that one expects $X$ to
contain a Zariski dense open subset of rational points, possibly defined
over some finite algebraic extension of $\Q$.
When $d\leq n$, the conjecture of Manin \cite{man}
predicts that one should have an asymptotic formula of the shape
$
N_U(B)=c_XB^{n+1-d}(1+o(1)),
$
as $B \rightarrow \infty$.  Here $U \subseteq X$ is the
open subset formed by deleting certain accumulating subvarieties
from $X$, and $c_X$ is a non-negative constant that has been given a
conjectural interpretation by Peyre \cite{peyre}.
Viewed in this light, our main result is most impressive in the case
$d=3$ and $n\geq 4$, in which setting one ought to be able to take $U=X$ in
Manin's conjecture.

The proof of Theorem \ref{6} is based upon 
an application of \cite[Theorem 4]{bhb}.
Broadly speaking this shows that every
point counted by $N_X(B)$ must lie on one of a small number of 
linear subspaces contained in $X$, each of which is defined over $\Q$.  
Thus one is naturally led to study the Fano variety
\beq\lab{fano-variety}
F_m(X)=\{\Lambda \in \Grass(m,n): \Lambda \subset X\},
\eeq
for $m\leq n-1$, 
where $\Grass(m,n)$ denotes the Grassmannian parametrising
$m$-dimensional linear subspaces $\Lambda\subset \bfP^n$. 
Perhaps the most basic example is provided by the case $m=1$ and
$d=n=3$, for which it is well-known that $F_1(X)$ has dimension $0$
and degree $27$.  The specific facts that we shall need are collected together in \S \ref{geom}. It
turns out that we have good control over the possible dimension of
$F_m(X)$ when $m=1$ or $m\geq (n-1)/2$, the latter fact being made
available to us by Professor Starr, in the appendix.

We end this introduction by summarising the contents of this paper.
The following section is concerned with detailing a number of basic estimates that
will be crucial to the proof of Theorem \ref{6}.   In particular we
shall need information about the growth rate of rational points on arbitrary projective varieties.
In \S \ref{geom} we shall collect together some facts about the geometry of non-singular hypersurfaces, and the
possible linear spaces that they contain.  An overview of the proof of Theorem
\ref{6} will be given in \S \ref{campaign}, before 
being carried out in full within \S\S \ref{line}--\ref{3-plane}.

\begin{ack}
While working on this paper, the first
author was supported at Oxford University by
EPSRC grant number GR/R93155/01.
The authors are grateful to the anonymous
referee, for his careful reading of the manuscript and numerous
pertinent suggestions. These led, in particular, to the discovery of a
significant error in the original proof of Lemma \ref{deg-cone}. 
\end{ack}

\section{Preliminary estimates}\lab{prelim}

In this section we collect together some of the basic estimates that
we shall need during the course of our work.  We begin with a rather
easy result from linear programming, whose proof we include for the
sake of completeness.

\begin{lem}\lab{prog}
Let $H \geq 1$, and let $a,b,c \geq 0$.  Then we have
$$
\max_{A,B,C} A^aB^bC^c \leq \max\{H^{(a+b+c)/3}, H^{(b+c)/2},H^c\},
$$
where the maximum on the left hand side 
is over all real numbers $A,B,C$ for which $1 \leq A
\leq B \leq C$ and $ABC \leq H$.
\end{lem}

\begin{proof}
Let $M=\max A^aB^bC^c$ in the statement of Lemma \ref{prog}.  Then on
writing $A=R$, $B=RS$ and $C=RST$, we see that  
$$
M= \max_{\colt{R,S,T\geq 1}{R^3S^2T\leq H}} 
R^{a+b+c}S^{b+c}T^c.
$$
Suppose first that $a+b \leq 2c$ and $b
\leq c$.  Then it follows that
$$
M\leq R^{3c}S^{2c}T^c \leq H^c.
$$
Next if $a+b >2c$ or $b > c$, then we substitute $T \leq
HR^{-3}S^{-2}$ and deduce that
$$
M \leq H^c \max_{\colt{R,S\geq 1}{R^3S^2\leq H}} R^{a+b-2c}S^{b-c}= M',
$$
say. Now it is easy to see that the maximum in the definition of $M'$ is
achieved at $S=1$ (resp. at $R=1$) if $a+b >2c$ and $b \leq c$ (resp. if
$a+b \leq 2c$ and $b > c$).  Thus 
$$
M' \leq \max\{H^{(a+b+c)/3}, H^{(b+c)/2}\}
$$
in either of these two case.  Finally if $a+b> 2c$ and $b > c$ then we
substitute $S \leq H^{1/2}R^{-3/2}$ and deduce that
$$
M' \leq H^{(b+c)/2} \max_{1\leq R^3\leq H} R^{a-b/2-c/2}\leq
\left\{
\begin{array}{ll}
H^{(b+c)/2}, & 2a \leq b+c,\\
H^{(a+b+c)/3}, & \mbox{otherwise}.
\end{array}
\right.
$$
This completes the proof of Lemma \ref{prog}.
\end{proof}

We shall also need some facts about 
the density of rational points on arbitrary locally closed subsets
$V\subset \bfP^N$.  We henceforth write $V(\Q)=V\cap \bfP^N(\Q)$ for
the set of rational points on $V$, and recall the definition 
(\ref{height}) of the projective height function
$H: \bfP^N(\Q)\rightarrow \R_{>0}$, given 
$x=[\x] \in \bfP^N(\Q)$ such that
$\x=(x_0,\ldots,x_N)\in \Z^{N+1}$ and $\hcf(x_0,\ldots,x_N)=1$.  
For any locally closed subset $V \subset \bfP^{N}$ and any
$B \geq 1$, we define the counting function
\begin{equation}\lab{flu}
N_V(B)=\#\{x \in V(\Q): ~H(x) \leq B\}.
\end{equation}
This coincides with our definition of $N_X(B)$ for a hypersurface
$X\subset\bfP^n$. 
When $V$ is a subvariety of $\bfP^N$ we shall always assume that 
it is defined over $\ov{\Q}$. 
Furthermore we shall henceforth refer to such a variety as being
integral if it is geometrically integral.
We then have the following ``trivial'' 
estimate, which is established in \cite[Theorem 1]{bhb}.

\begin{lem}\lab{triv}
Let $V \subset \bfP^{N}$ be a variety of degree $d$ and
dimension $m$.   Then we have
$$
N_V(B) =O_{d,N}(B^{m+1}).
$$
\end{lem}

It is easy to see that Lemma \ref{triv} is best possible when $V$
contains a linear subspace of dimension $m$ that is defined over $\Q$.
On the assumption that $V$ is integral and has degree $d \geq
2$, we can do somewhat better than Lemma \ref{triv}.  
The following result is extracted from the the introduction to
\cite{bhbs}.

\begin{lem}\lab{pila'}
Let $\ve>0$ and suppose that $V \subset \bfP^N$ is an integral 
variety of degree $d\geq 2$ and dimension $m$.  
Then we have
$$
N_V(B)\ll_{d,\ve,N}\left\{
\begin{array}{ll}
B^{m+1/4+\ve}, & \mbox{if $m \geq 4$ and $3\leq d \leq 5$},\\
B^{m+\ve}, & \mbox{otherwise}.\\
\end{array}
\right.
$$
\end{lem}

It will be clear to the reader that when
$m \geq 4$ and $3\leq d \leq 5$, the main result in \cite{bhbs}
actually allows one to take a sharper exponent in the statement of
Lemma \ref{pila'}.  In fact the exponent $m+\delta_d$ is acceptable
for any
$$
\delta_d>
\left\{ \begin{array}{ll}
5/(3\sqrt{3})-3/4, &d=3,\\
1/12, & d\geq 4.
\end{array}
\right.
$$
However it turns out that
the estimate provided above is sufficient for the purposes of
Theorem \ref{6}.

For non-negative integers $m\le N$, let $\Grass(m,N)$ denote the Grassmannian
which parametrises $m$-planes contained in $\bfP^{N}$. 
It is well-known that $\Grass(m,N)$ can be embedded in $\bfP^\nu$ via
the Pl\"ucker embedding, where
$$
\nu=\binom{N+1}{m+1}-1,
$$
and that $\Grass(m,N)$ has dimension 
$(m+1)(N-m).$  If $M \in \Grass(m,N)(\Q)=\Grass(m,N)\cap
\bfP^\nu(\Q)$, we define the height $H(M)$ of
$M$ to be the standard multiplicative height of its coordinates in
$\Grass(m,N)$, under the Pl\"ucker embedding.   
The following result is well-known (see 
\cite[\S 2.4]{3fold1}, for example),
and refines Lemma~\ref{triv} in the case of
linear varieties.

\begin{lem}\lab{triv-lin}
Let $M \in \Grass(m,N)$.  Then we have
$$
N_M(B)  \ll_{N} B^{m}+\frac{B^{m+1}}{H(M)}.
$$
Moreover, if $M$ contains $m+1$ linearly independent rational
points of height at most $B$, then $M$ is defined over $\Q$ and 
$$
\frac{B^{m+1}}{H(M)}\ll_N N_M(B)  \ll_{N} \frac{B^{m+1}}{H(M)}.
$$
\end{lem}

The following key result  allows us to tackle Theorem \ref{6} by restricting 
attention to the linear spaces that are contained in the
hypersurface $F=0$.

\begin{lem}\lab{main1}
Let $\ve>0$ and suppose that $X\subset \bfP^{n}$ is an
integral hypersurface of degree $d$.
Then there exist linear spaces $M_1, \ldots, M_J \subseteq X$ defined
over $\Q$, with $J =O_{d,\ve,n}(B^{n-1+\ve})$,
such that $0\leq \dim M_j\leq n-1$ for $1\leq j \leq J$ and 
$$
N_X(B) \leq \sum_{j=1}^J N_{M_j}(B).
$$
Moreover whenever $\dim M_j \geq 1$, for any $1 \leq j \leq J,$ we have
$$
H(M_j)=O_{d,\ve,n}(B^{1+\ve}).
$$
\end{lem}

\begin{proof}
Since $X$ is integral there exists an 
absolutely irreducible form $F \in \ov\Q[X_0,\ldots,X_n]$ of degree
$d$, such that $X$ is given by the equation $F=0$.
We first assume that $F$ is not proportional to a
form defined over $\Q$ and let $F^\sigma$ be the conjugate of $F$ for
any non-trivial $\sigma \in \rom{Gal}(\ov{\Q}/\Q)$.  Then
clearly 
$
N_X(B) =N_{X\cap X^{\sigma}}(B),
$ 
where $X^\sigma$ is the
hypersurface $F^\sigma=0.$  Since $X\cap X^{\sigma} \subset
\bfP^{n}$ is a variety of dimension at most $n-2$, we have 
$$
N_X(B)\ll_{d,n}B^{n-1},
$$  
by Lemma \ref{triv}.  Thus we
may take $M_1,\ldots,M_J$ to be the collection of zero dimensional 
linear subspaces that
correspond to precisely these points, in the statement of Lemma
\ref{main1}.  

Suppose now that $F \in \Z[X_0,\ldots,X_n]$,
and let $\|F\|$ denote the maximum modulus of the
coefficients of $F$.  We claim that there exists a constant $c_{d,n}$
depending only on $d$ and $n$,     
such that if $\log \|F\| >
c_{d,n} \log B$, then there is a hypersurface $Y \subset \bfP^{n}$ of degree $d$, different from 
$X$, such that $N_X(B) =N_{X\cap Y}(B)$.  But this is a direct
consequence of \cite[Lemma 3]{bhb}.   Thus 
the case in which $F \in \Z[X_0,\ldots,X_n]$ has 
$\log \|F\| > c_{d,n} \log B$ is also
satisfactory for Lemma \ref{main1}, by Lemma \ref{triv}.

Finally we suppose that $F \in \Z[X_0,\ldots,X_n]$ and that $\log \|F\|
\ll_{d,n} \log B$.  But then a direct application of \cite[Theorems 4
and 5]{bhb} yields  the result.
The reader should note that \cite{bhb} works with forms in
$\Z[X_1,\ldots,X_n]$, rather than in $\Z[X_0,\ldots,X_n]$, so that the
two values of the parameter $n$ do not correspond.
\end{proof}

Let $m \in \N$ and let $V \subset \bfP^N$ be an integral variety of degree $d$.
Our final result in this section provides a crude upper bound for the
dimension of the Fano variety
$$
F_m(V)=\{\Lambda \in \Grass(m,N): \Lambda \subset V\},
$$
which parametrises the $m$-planes contained in $V$.  In the next 
section we shall
see how much more can be said when $V$ is a
non-singular hypersurface.  The following estimate may certainly be
extracted from the work of Segre \cite{segre}, although we have provided
our own proof for the sake of completeness.

\begin{lem}\lab{gen-fano}
We have
$$
\dim F_1(V) \leq 
\left\{
\begin{array}{ll}
2\dim V -2, & d=1,\\
2\dim V -3, & d \geq 2.
\end{array}
\right.
$$
\end{lem}
\begin{proof}
Let $\delta=\dim V$.  In order to establish Lemma \ref{gen-fano} 
we note that the case in which
$V$ is isomorphic to $\bfP^\delta$ is easy, since then
$F_1(V)\cong\Grass(1,\delta)$. Assuming therefore that
$d \geq 2$, we employ a routine incidence correspondence argument.
Let $Z$ be an integral component of $F_1(V)$, and let
\[\Sigma=\{(v,L) \in V \times Z:\; v \in L\}.\]
Then consideration of the projection onto the second factor shows that
$\dim \Sigma = \dim Z +1$.  Now let $V_0\subseteq V$ be the union of the 
lines in $Z$, and let $v$ be a generic point of $V_0$.  The projection 
$\Sigma \rightarrow V$ then shows that $\dim \Sigma=\dim V_0+\dim Z_v$,
where $Z_v=\{L \in Z: v \in L\}$.  Thus
\[\dim Z=\dim \Sigma-1=\dim V_0-1+\dim Z_v\le\dim V-1+\dim Z_v.\]
Now any line $L\in Z_v$ must
also lie in the tangent space $\mathbb{T}_v(V_0)$, so that
$Z_v\subseteq W_v$, where 
\[W_v=\{L\in \Grass(1,N): v\in L\subseteq\mathbb{T}_v(V_0)\}.\]
Since $v$ is generic on $V_0$ it is non-singular, so that
$W_v$ is a linear space of dimension $\dim V_0-1$.

We proceed to consider two cases.  If $V_0$ is a linear space then it must be
a proper subvariety of $V$, since $d\ge 2$.  In this case 
\[\dim Z_v\le\dim W_v=\dim V_0-1\le\dim V-2,\]
and the required result follows.  On the other hand if $V_0$ is not
linear, then $Z_v$ must be a proper subvariety of $W_v$, and
\[\dim Z_v\le\dim W_v-1=\dim V_0-2\le\dim V-2.\]
Again this suffices for the lemma.
\end{proof}

\section{Geometry of non-singular hypersurfaces}\lab{geom}

For any $d \geq 3$ and $n \geq 4$, let $X \subset \bfP^{n}$ 
be a non-singular projective hypersurface of
degree $d$. 
The aim of this section is to discuss the geometry of such
hypersurfaces, and in particular the possible $m$-planes in 
$\bfP^{n}$ that are contained in them.
Such $m$-planes are parametrised by the Fano
variety $F_m(X)$, given by (\ref{fano-variety}).
We begin by collecting
together some preliminary facts about the degree and dimension of $F_m(X)$,
for various values of $m \in \N$.

\begin{lem}\lab{fano}~
\begin{enumerate}
\item[(i)] $F_m(X)$ has degree  $O_{d,n}(1)$ when it is non-empty.
\item[(ii)] $F_m(X)$ is empty for $m> (n-1)/2$.
\item[(iii)]$F_m(X)$ is finite if $m= (n-1)/2$.
\item[(iv)] $F_1(X)$ has dimension $2n-3-d$ when $d\leq \min\{6,n\}$.
\item[(v)]  $F_2(X)$ has dimension $3n-16$ for $d=3$ and $n \geq 6$.
\end{enumerate}
\end{lem}

\begin{proof}
The first part of Lemma \ref{fano} is well known, and follows 
from using the defining equation for $X$ to write down the explicit 
equations for $F_m(X)$.
In the appendix
Professor Starr has provided proofs of (ii) and (iii).  As indicated there it
is very easy to establish part (ii), whereas the boundary case
in part (iii) requires subtler methods.
Part (iv) is the principal result in recent work of  Beheshti 
\cite{zav}, and part (v) follows from work of Izadi \cite[Prop. 3.4]{izadi}.
\end{proof}

In the case $m=1$ of lines, it is interesting to remark that a
standard incidence correspondence argument (see Harris 
\cite[\S 12.5]{harris}, for example) reveals that 
$\dim F_1(X)=2n-3-d$ for a generic
hypersurface $X \subset \bfP^{n}$ of degree $d$. 
Debarre and de Jong have  conjectured this to be the true
dimension of $F_1(X)$ whenever $d \leq n$. 
When $m=2$ the situation seems to be less well understood.  While the
dimension of $F_2(X)$ is $3n-6-(d+1)(d+2)/2$ for a generic
hypersurface $X \subset \bfP^{n}$ of degree $d$, which tallies with
part (v) of Lemma \ref{fano}, there exist examples showing  that this
is not always the true dimension when $d\leq n$.  Perhaps the simplest
example is provided by the Fermat cubic in $\bfP^5$ which contains a
finite number of planes.  Nonetheless, we shall see below in Lemma
\ref{fano'} that it is possible to prove non-trivial upper bounds for
the dimension of $F_2(X)$ when $d=4$ and $X$ is covered by planes.

It is convenient at this point to raise a question concerning the
possible dimension of $F_m(X)$ in the 
sub-boundary case $m=\lfloor(n-1)/2\rfloor$, where $\lfloor\alpha\rfloor$ denotes the integer
part of $\alpha \in \R$.  During the course of our work we 
have been led to formulate the
following conjecture, the resolution of which would simplify
the proof of Theorem \ref{6} considerably.

\begin{con}
\lab{hope}
For any $n\geq 5$, let $X \subset \bfP^{n}$ be a non-singular 
hypersurface of degree
$d\geq 3$, and let $m=\lfloor(n-1)/2\rfloor$.
Then $X$ is not a union of $m$-planes.\footnote{Since this paper was
submitted for publication, the conjecture has been proved when $d\geq
  4$ by Roya Beheshti \cite{zav'}.}
\end{con}

Conjecture \ref{hope} is already included in part (iii) of
Lemma~\ref{fano} when $n$ is odd, since $X$ contains at most finitely
many linear subspaces of dimension $m$ in this case.  
Further evidence is provided by part (v) of the same result.  Indeed when $d=3$ and
$n=6$ we see that $F_2(X)$ has dimension $2$, so that
$
\bigcup_{P \in F_2(X)}P
$
has dimension at most $4$ and must be a proper subvariety of $X$.
It is clear that Conjecture \ref{hope} can be false when $n\leq 4$, 
since for example a non-singular cubic threefold is covered by its lines.

For any $m \in \N$, and any subvariety  $Z \subseteq F_m(X)$, we
shall henceforth write
\beq\lab{C}
D(Z)=\bigcup_{\Lambda \in Z}\Lambda,
\eeq
to denote the union of $m$-planes swept out by $Z$.
It follows from part (ii) of Lemma \ref{fano} that $D(Z)$ 
is empty for $m>(n-1)/2$.  Clearly $D(Z) \subseteq X$ and it is not
hard to see that the degree of $D(Z)$ is bounded in terms of
$d,n$ and the degree of $Z$.  We shall use these facts
without further comment throughout this paper.

We proceed by discussing various cones of $m$-planes contained in
$X$.  For any subvariety $Z \subseteq F_m(X)$, and any $x \in X$, we set
\beq\lab{fibred}
Z_x=\{\Lambda \in Z: x \in \Lambda\}.
\eeq
We shall write
$$
C_x(Z)=D(Z_x)
$$
for the corresponding cone of $m$-planes through $x$.
In particular, since the degree of $Z_x$
is bounded in terms of $d,n$ and the degree of $Z$, it follows from the
previous paragraph that the degree of $C_x(Z)$ is also bounded in terms of
$d,n$ and the degree of $Z$.
We now observe that if a hypersurface of degree at least $2$ is a
cone, then its vertex is a singular point. Thus, since $X$ is
non-singular and $C_x(Z) \subseteq X$, we must have  
\beq\lab{condition}
\dim C_x(Z) \leq n-2
\eeq
for any $x \in X$ and any subvariety $Z \subseteq F_m(X)$.  

Suppose now that $m = 2$, and let $L \in F_1(X)$.  Then analogously
to (\ref{fibred}), we define $Z_L=\{P \in Z: L \subset
P\}$ and $C_L(Z)=D(Z_L)$ for any closed subset $Z \subseteq
F_2(X)$. The following result corresponds to (\ref{condition}).

\begin{lem}\lab{condition'}
Let $d \geq 3$ and let $n \geq 5$. 
Then we have 
$$
\dim C_L(Z) \leq n-3,
$$
for any $L \in F_1(X)$ and any subvariety $Z\subseteq F_2(X)$. 
\end{lem}

\begin{proof}
For fixed $L \in F_1(X)$ it will clearly suffice to establish the
upper bound
\beq\lab{n-5}
\dim Z_L \leq n-5,
\eeq
under the assumption that $X \subset \bfP^n$ is a non-singular
hypersurface of degree $d \geq 3$ and dimension $n-1\geq 4$. 
Pick any distinct points $x,y \in L$ and let $P\in Z_L$.  
Then the system of inclusions 
$x\in L\subset P\subset X$ implies that $L\subset P\subset
\mathbb{T}_x(X)$, where $\mathbb{T}_x(X)\cong \bfP^{n-1}$ is the
tangent hyperplane to $X$ at $x$.  Similarly we have 
$L\subset P\subset \mathbb{T}_y(X)$.
Hence it follows that
\begin{align*}
Z_L&=\{P\in Z: L\subset P \subset \mathbb{T}_x(X)\cap\mathbb{T}_y(X)\}\\
&=Z\cap\{P\in F_2(\mathbb{T}_x(X)\cap\mathbb{T}_y(X)): L\subset P\}\\
&=Z \cap G_L,
\end{align*}
say.  

We claim that it is possible to choose $x$ and $y$ in $L$ so 
that the tangent spaces $\mathbb{T}_x(X)$ and $\mathbb{T}_y(X)$ are
different. Suppose the hypersurface $X$ is defined by the form
$F(X_0,\ldots,X_n)$, and let $x=[\ma{x}]$ and $y=[\ma{y}]$ be distinct
points on $L$.  If the tangent space is the same for all points on $L$
then there is a fixed vector $\ma{v}$, say, such that
$\nabla F(\lambda\ma{x}+\mu\ma{y})$ is always a scalar multiple of
$\ma{v}$.  In fact we must have
$\nabla F(\lambda\ma{x}+\mu\ma{y})=h(\lambda,\mu)\ma{v}$ 
for a certain binary form $h$ of degree $d-1$.  Thus there is a pair
$(\lambda,\mu)\not=(0,0)$ for which $h(\lambda,\mu)=0$.  We therefore
obtain a singular point on the variety $X$.  This contradiction
establishes our claim.

Now, with an appropriate choice of $x$ and $y$, we see that
\[G_L \cong \{P\in \Grass(2,n-2): L\subset P\}
\cong\bfP^{n-4}.  \]
In order to complete the proof of (\ref{n-5}) it
plainly suffices to show
that $G_L$ is not contained in $Z$. Arguing by contradiction we suppose that
$G_L\subseteq Z$ and form the union of planes $D(G_L)$.
But then $D(G_L) \subset X$ is a linear algebraic variety of dimension $n-2$,
which is impossible by part (ii) of Lemma \ref{fano}.
\end{proof}

We now come to the most important 
lemma in this section --- a kind of stratification
result that will form the backbone of our proof of Theorem \ref{6}.
Let $m \in \N$, let $\Phi \subseteq F_m(X)$ be an integral
component and let $Y=D(\Phi) \subseteq X$ so that $Y$ is also
integral. We proceed by considering the incidence correspondence
\beq\lab{I}
I=\{(y,\Lambda) \in Y \times \Phi: y \in \Lambda\}.
\eeq
Then the projection onto the first factor is surjective, and
by projecting onto the second factor we see that
$I$ has  dimension $\dim \Phi+m$.  Thus it follows that
\beq\lab{gen-dim}
\dim \Phi_{y}= \dim \Phi-\dim Y+m
\eeq
for generic $y \in Y$, in the notation of (\ref{fibred}).
The following result shows how the dimension of
$\Phi_y$ varies for different choices of $y \in Y$. We henceforth employ 
the convention that the empty set is the only algebraic set
with negative dimension.

\begin{lem}\lab{strat}
Let $\Phi \subseteq F_m(X)$ be integral and let $Y=D(\Phi)
\subseteq X$ have degree $e$. 
Then there exists a stratification of subvarieties
$$
Y=Z_0(\Phi) \supseteq Z_1(\Phi) \supseteq Z_2(\Phi) \supseteq \cdots,
$$
such that the following holds.
\begin{enumerate}
\item[(i)]
For $i \geq 1$ we have $\deg Z_i(\Phi)=O_{e,n}(1)$ and
$$
\dim Z_i(\Phi) \leq \dim Y-1-i.
$$
\item[(ii)]
For $i \geq 0$ and any
$y \in Z_i(\Phi) \setminus Z_{i+1}(\Phi)$ we have
$$
\dim \Phi_y = \dim \Phi-\dim Y+m+i.
$$
\end{enumerate}
\end{lem}

\begin{proof}
Throughout this proof we shall write $\del$ for the dimension of
$\Phi$, and $\ve$ for the dimension of $Y$.
Recall the definition (\ref{I}) of the incidence correspondence $I$.
We have already seen that $I$ has  dimension $\del+m$, and that
(\ref{gen-dim}) holds for generic $y \in Y$. 
In order to make this more precise we define 
$Y_k$ to be the set of $y \in Y$ for which
$\dim \Phi_y \geq k$, for any non-negative integer $k$.
Then $Y_k$ is a closed subset of
$Y$ and has degree $O_{e,n}(1)$.  The first statement here follows
from the upper semicontinuity of the
dimension of the fibres of a morphism (see \cite[Corollary
11.13]{harris}, for example), while the latter fact can be extracted
from an analysis of the proof of \cite[Theorem 11.12]{harris}.
On noting that $\dim \Phi_y \leq \delta$ for any $y \in Y$, 
it plainly follows that 
$$
\emptyset = Y_{\del+1} \subseteq Y_{\del} \subseteq \cdots \subseteq Y_0=Y.
$$
Moreover $Y_{\del -\ve+m+i}$ is a
proper subvariety of $Y$ for $i \geq 1$.

We shall take 
$$
Z_0(\Phi)=Y, \quad Z_i(\Phi)=Y_{\del-\ve+m+i}, 
$$
for $i \geq 1$, in the statement of Lemma \ref{strat}. Then part (ii)
of the lemma is immediate, and it
remains to provide an upper bound for the dimension of $Z_i(\Phi)$ when
$i \geq 1$. For this we consider the incidence correspondence 
$$
I_{i}= \{(y,\Lambda) \in Z_{i}(\Phi) \times \Phi: y \in \Lambda\}.
$$
Since (\ref{gen-dim}) holds for generic $y\in Y$ it follows that
$Z_i(\Phi)$ is a proper subvariety of $Y$ for $i\geq 1$.  Hence the
generic $m$-plane $\Lambda\in \Phi$ cannot lie completely in $Z_i(\Phi)$.
Now let $\pi_2$ be the projection from $I_i$ to $\Phi$.  If $\pi_2$ 
is onto, then the
generic fibre $\{y\in Z_i(\Phi)\cap\Lambda\}$ has dimension at most
$m-1$, in which case $\dim I_i \leq \del+m-1.$
On the other hand, if $\pi_2$ is not onto, then $\dim
\pi_2(I_i)\le\del-1$ and we again deduce that 
\[\dim I_i \leq \del+m-1. \]

Turning to the projection to the first factor, we see that it 
is onto, since $Z_i(\Phi)\subseteq Y=D(\Phi)$, and hence
\[\dim I_i\ge \dim Z_{i}(\Phi) + \del-\ve+m+i.\]
Thus it follows that
$$
\dim Z_{i}(\Phi) \leq \ve-1-i
$$
for $i \geq 1$, as required.
This completes the proof of Lemma \ref{strat}.
\end{proof}

In the special case $m=1$ and $Z=F_1(X)$, we can actually write down
the equations defining the cone $C_x(Z)$ for any $x \in X$.
Let $C_x^1=C_x(F_1(X))$ for a
point $x \in X$, which we assume without loss of generality is given by
$x=[1,0,\ldots,0]$.  Then, after a further linear change of variables, 
$X$ takes the shape
\beq\lab{X}
X_0^{d-1}X_1+ X_0^{d-2}F_2(X_1,\ldots,X_n)+\cdots+F_{d}(X_1,\ldots,X_n)=0,
\eeq
for forms $F_i$ of degree $i$, for $2\leq i \leq d$.
If $x\in L \in F_1(X)$, then $L$ must be contained in 
the tangent hyperplane $\mathbb{T}_x(X)$, which is given by the equation
$X_1=0$.  It follows that any point in $C_x^1$ must
be of the form $[a,0,b_2,\ldots,b_n]$, and since 
$C_x^1\subset X$ the polynomial
$$
a^{d-2}F_2(0,b_2,\ldots,b_n)+\cdots+F_{d}(0,b_2,\ldots,b_n)
$$
must vanish identically in $a$. We therefore deduce that
\beq\lab{eqn}
C_x^1=\{[a,0,\ma{b}] \in \bfP^n:  F_2(0,\ma{b})=\cdots=F_d(0,\ma{b})=0\},
\eeq
where $\ma{b}=(b_2,\ldots,b_n)$.  We see from (\ref{eqn}) that $\dim
C^1_x\ge n-d$.  Moreover it is plain that
$C_x^m=C_x(F_m(X))\subseteq C_x^1$ for any $m \in \N$.

When $d=3$ and $x\in X$ is generic we can be even more precise about the
cone $C_x^1$ as soon as $n$ is
large enough.  Given a non-singular cubic hypersurface $X \subset
\bfP^n$, Barth and Van de Ven \cite[Prop. 5]{barth} have shown that 
$F_1(X)$ is a non-singular simply connected variety
of dimension $2n-6$, for $n \geq 6$.  
The fact that $F_1(X)$ is a non-singular variety of dimension $2n-6$,
for $n \geq 6$, is already present in the work of Clemens and Griffiths \cite{clem},
where it is also shown that $X$ is a union of lines.  We may conclude that
$F_1(X)$ is integral for $n \geq 6$, and we proceed to establish the following result.

\begin{lem}\lab{deg-cone}
Let $d=3$ and let $n \geq 7$.  Then for generic $x \in X$ the 
cone $C_x^1$ is integral and has degree $6$.
\end{lem}

\begin{proof}
We begin by showing that $C_x^1$ has dimension $n-3$ for generic $x \in
X$. This is an easy consequence of (\ref{gen-dim}), since we
can take $\Phi=F_1(X)$ and $Y=X$. 
Thus for a generic point $x\in X$ one has
\[
\dim \Phi_x=\dim F_1(X)-\dim X+1=n-4,
\]
by part (iv) of Lemma \ref{fano}.  It follows that $\dim
C^1_x=1+\dim\Phi_x=n-3$, as claimed.
We shall also make use of the fact for a non-singular cubic
hypersurface $X \subset \bfP^n$, the Hessian $H$ of $X$ does not
contain $X$ as a subvariety.  This is established by Hooley
\cite[Lemma 1]{hooley}, for example, and implies in particular that
$H\cap X$ is a proper subvariety of $X$.

To find the degree of $C^1_x$ for generic $x \in X$ it is convenient
to choose coordinates as in (\ref{X}) and (\ref{eqn}).  Thus we may
assume without loss of generality that $x=[1,0,\ldots,0]$, and $X$ is
defined by the non-singular cubic form
\beq\lab{vault}
F(X_0,\ldots,X_n)=X_0^2X_1 + X_0Q(X_1,\ldots,X_n)+C(X_1,\ldots,X_n),
\eeq
for some quadratic and cubic forms $Q$ and $C$, respectively, and furthermore
$$
C_x^1=\{[a,0,\ma{b}] \in \bfP^n:  
Q(0,\ma{b})=C(0,\ma{b})=0\}.
$$ 
Since $C_x^1$ has dimension $n-3$, it follows that neither 
$Q(0,\ma{b})$ nor $C(0,\ma{b})$
vanish identically, and that $C_x^1$ is a complete intersection of
pure dimension $n-3$.  In fact $Q(0,\ma{b})$ must be a non-singular
quadratic form, since the Hessian of
$F$ at $x$ is equal to $-4 \det(Q(0,\ma{b}))$, and $x\in X$ is generic.
We shall think of $C_x^1$ as lying in
$\bfP^{n-1}$, by identifying $[a,0,\ma{b}]$ with $[a,\ma{b}]$.

In order to complete the 
proof of Lemma \ref{deg-cone}
it suffices to show that $C_x^1$ is
reduced and irreducible, by B\'ezout's theorem.
Suppose for a contradiction that $C_x^1=Y_1\cup Y_2$, for components
$Y_1, Y_2$ of dimension $n-3$ in $\bfP^{n-1}$.  Then any point in the
intersection $Y_1 \cap Y_2$ produces a singular point on $C_x^1$, so
that the singular locus of $C^1_x$ has dimension at least $n-5$.
Now write $\nabla'=(\partial/\partial X_2,\ldots,\partial/\partial
X_n)$, and let $V$ denote the set of $b=[\ma{b}] \in \bfP^{n-2}$ such
that $Q(0,\ma{b})=C(0,\ma{b})=0$ and $\nabla' Q(0,\ma{b})$ is
proportional to $\nabla' C(0,\ma{b})$.  Then $V$ is an 
algebraic subvariety of $\bfP^{n-2}$, since the latter condition is
defined by the vanishing of various $2\times 2$ determinants.
Moreover it is clear that $\dim V\geq n-6$.

For $1\leq i\leq n$, let $C_i$ (resp. $Q_i$) denote the partial derivative 
$\partial Q/\partial X_i$ (resp. $\partial C/\partial X_i$).
We proceed to define the map $\pi:V \rightarrow \bfP^{n-1}$,  via
$$
\pi: [\ma{b}]\longmapsto [C_i(0,\ma{b}),Q_i(0,\ma{b})b_2,\ldots,Q_i(0,\ma{b})b_n],
$$
whenever $C_i(0,\ma{b}), Q_i(0,\ma{b})$ do not both vanish.
It is not hard to check that $\pi$ is well-defined, since
$Q(0,\ma{b})$ is non-singular.  We claim that $\pi(V)$ has dimension
at least $n-6$, for which it is clearly enough to show that $\pi$ is
generically injective.  But if $[u,\ma{v}]$ is a generic
point in the image $\pi(V)$, then we cannot have $\ma{v}=\ma{0}$ since 
$Q(0,\ma{b})$ is non-singular. Thus $[u,\ma{v}]$ determines
$[\ma{b}]$, and so $\dim \pi(V)\geq n-6$, as required.  
Now any point $[u,\ma{v}]\in \pi(V)$ must satisfy
$$
Q(0,\ma{v})=0, \quad u \nabla' Q(0,\ma{v})= \nabla' C(0,\ma{v}),
$$
and it follows that the set $W$ of such $[u,\ma{v}]\in \bfP^{n-1}$
has dimension at least $n-6$.  
Finally we may conclude that the set of $[u,\ma{v}]\in W$ for which 
$$
u^2+u \frac{\partial Q}{\partial X_1}(0,\ma{v})+
\frac{\partial C}{\partial X_1}(0,\ma{v})=0,
$$
has dimension at least $n-7.$ 
Since $n \geq 7$, we produce at least one point $(u,0,\ma{v})$ at 
which the form (\ref{vault}) is singular. 
This contradiction completes the proof of Lemma \ref{deg-cone}.
\end{proof}

We now turn to the special case of planes contained in
non-singular quartic hypersurfaces, with a view to proving the result
alluded to in the paragraph after Lemma \ref{fano}.
With this in mind we have the following result.

\begin{lem}\lab{fano'}
Let $d=4$ and let $n \geq 6$.  Then for 
any integral component $\Phi \subseteq F_2(X)$ such
that $X=D(\Phi)$, we have
$$
\dim \Phi \leq \left\{
\begin{array}{ll}
3, & n=6,\\
3n-16, & n \geq 7.
\end{array}
\right.
$$
\end{lem}

\begin{proof}
As above we let $C_x^1=C_x(F_1(X))$ denote the union of all lines
contained in $X$ that pass through a point $x \in X$.  For generic $x \in
X$ we claim that 
\beq\lab{needed'}
\dim C_x^1 = n-4,
\eeq 
and that if $n \geq 7$ the cone $C_x^1$ does not contain 
any linear space of dimension $n-4$.
The latter claim follows from the fact $X$ contains at most finitely
many $(n-4)$-planes by parts (ii) and (iii) of Lemma \ref{fano}, so
that a generic point of $X$ cannot lie on such a subvariety.  To
see the first claim, one notes that if $\Psi \subseteq F_1(X)$ is any integral
component such that $D(\Psi)$ is a proper subvariety of
$X$, then $C_x(\Psi)$ is empty for generic $x \in X$.
Alternatively, if $\Psi \subseteq F_1(X)$ is an integral
component such that $D(\Psi)=X$ then 
one combines (\ref{gen-dim}) with part (iv) of Lemma \ref{fano},
just as in the proof of Lemma \ref{deg-cone}, to deduce (\ref{needed'}).

We proceed by considering the cone of planes $C_x(\Phi)=D(\Phi_x)$,
for any $x \in X$.   In particular $C_x(\Phi) \subseteq C_x^1$.
Moreover if $H \in {\bfP^{n}}^*$ is a generic
hyperplane, then $H$ does not contain $x$ and it follows
that $H \cap P$ is a line for every $P \in \Phi_x$. 
Thus there is a bijection between planes parametrised by $\Phi_x$ and
lines contained in $H \cap C_x(\Phi)$, whence
$$
\dim \Phi_x= \dim F_1(H \cap C_x(\Phi))
$$
for any $x \in X$.  Now let $x \in X$ be generic.  Then on combining
our observation that $C_x(\Phi) \subseteq C_x^1$ with
(\ref{needed'}), we conclude that
\beq\lab{needed''}
\dim \Phi_x\leq  \dim F_1(H \cap C_x^1),
\eeq
where $H \cap C_x^1 \subset \bfP^{n}$ is a variety of dimension
$n-5$.  Moreover, since $H$ is
generic, the only way that $H \cap C_x^1$ can contain an $(n-5)$-plane
is if $C_x^1$ contains an $(n-4)$-plane.  We have already seen that
this is impossible when $x \in X$ is generic and $n \geq 7$.  
On applying Lemma \ref{gen-fano} to (\ref{needed''}) we therefore deduce
that
$$
\dim \Phi_x \leq 2(n-5)-2 = 0
$$
if $n=6$, whereas 
$$
\dim \Phi_x \leq 2(n-5)-3 = 2n-13
$$
if $n \geq 7$. To complete the proof of Lemma \ref{fano'} it now
suffices to apply (\ref{gen-dim}) with $Y=X$ and $m=2$.
\end{proof}

It is likely that the upper bound in Lemma
\ref{fano'} is not best possible, and the problem of proving sharper
versions seems to be an interesting question in its own right.

The final result of this section pertains to the special case $m=3,
n=8$.   We shall use the notation $C_x^3=C_x(F_3(X))$ for $x \in X$.
We have already seen in (\ref{condition}) that $\dim C_x^3 \leq n-2$.
The following result investigates when the
dimension of $C_x^3$ is maximal.

\begin{lem}\lab{final}
Let $d=3$ or $4$ and let $n=8$. Suppose that $\dim C_x^{3}=6$ for
some $x \in X$.  Then we must have $d=4$ and $\dim G_x=3$, where
$$
G_x= \{T \in F_3(X): x\in T\}.
$$
\end{lem}

\begin{proof}
Let $x \in X$ be such that 
$\dim C_x^{3}=6$, and note that $C_x^3=D(G_x)$ in
the notation of (\ref{C}).  In particular it follows that $\dim G_x
\geq 3$.  Assuming that $x=[1,0,\ldots,0]$, and that $X$ takes the
shape (\ref{X}), it
follows from (\ref{condition}) and (\ref{eqn}) that 
$$
C_x^3\subseteq C_x^1=\{[a,0,\ma{b}] \in \bfP^8:  
F_2(0,\ma{b})=\cdots=F_d(0,\ma{b})=0\},
$$ 
where $\dim C_x^1 \leq 6$.  If we define the variety
$$
Y=\{[\ma{b}] \in \bfP^6:  F_2(0,\ma{b})=\cdots=F_d(0,\ma{b})=0\} 
\subset \bfP^6,
$$
then it follows that $Y$ must have dimension $5$, 
since $C_x^3$ has dimension $6$.  Hence there exists a non-constant form $H \in
\ov{\Q}[X_2,\ldots,X_8]$ such that  for $2 \leq i \leq d$ we have 
$H \mid F_i(0,X_2,\ldots,X_8)$.
Thus the equation for $X$  takes the shape 
$$
F(X_0,\ldots,X_8)=
X_1J(X_0,\ldots,X_8) +H(X_2,\ldots,X_8)K(X_1,X_2,\ldots,X_8)=0,
$$
for some further forms $J$ and $K$ such that $\deg J=d-1$ and
$\deg H+\deg K=d$.
This is clearly impossible unless $K$ is in
fact a constant, since otherwise we may find a common solution to the
system of equations $X_1=J=H=K=0$, and this will produce a singular
point on $X$.  Thus $F_2(0,\ma{b}),\ldots,F_{d-1}(0,\ma{b})$ must all
vanish identically.  Taking the constant $K$ to be $1$, we then see that
$X$ is given by
\[F(X_0,\ldots,X_8)=X_1J(X_0,\ldots,X_8) +H(X_2,\ldots,X_8)=0,\]
and that $Y$ is a degree $d$ hypersurface in $\bfP^6$, given by
$H(X_2,\ldots,X_8)=0$.

We proceed to show that $Y$ is non-singular.  In general we have
$$
\nabla F = \Big(X_1\frac{\partial J}{\partial X_0}, J+X_1\frac{\partial
J}{\partial X_1}, X_1\nabla' J+ \nabla' H\Big),
$$
where 
$$
\nabla'=(\partial/\partial X_2, \ldots,\partial/\partial X_8).
$$
Now, given any non-zero vector $\x=(x_2,\ldots,x_8)$ corresponding to
a singular point on $Y$, we will have $\nabla' H(\x)=\ma{0}$. From
this it follows from Euler's identity that $H(\x)=0$.  We can 
then find $y \in \ov\Q$ such that $\nabla F(y,0,\ma{x})=\ma{0}$.
This produces a singular point on $X$, which is
impossible, thereby establishing that
$Y \subset \bfP^6$ is a non-singular hypersurface of degree $d.$  
We shall write
\[Z=\{[0,0,\ma{b}] \in \bfP^8: [\ma{b}]\in Y\},\]
so that $C_x^3$ is a cone over $Z$.

Let $T \in G_x$ and define $\pi(T)=T \cap Z$.  
Then it is clear that $T \cap Z \subset Z$ must be a 2-plane for each $T
\in G_x$, so that we have a map $\pi: G_x \rightarrow F_2(Z)$.
Since $x \not\in Z$, for each given $P \in F_2(Z)$ there is a 
unique $3$-plane which contains $P$ and also passes through $x $.
It follows that $\pi$ is a bijection.
We may now deduce from Lemmas \ref{fano} and \ref{fano'} that 
$$
\dim G_x=\dim F_2(Z)=\dim F_2(Y)\leq
\left\{
\begin{array}{ll}
2, & d=3,\\
3, &d=4,
\end{array}
\right.
$$
since $Y\subset \bfP^6$ is a non-singular hypersurface of degree $d$
and $D(F_2(Y))=Y$.  In particular this is impossible when $d=3$
since we have already seen that $\dim G_x \geq 3$.
This completes the proof of Lemma \ref{final}.
\end{proof}

\section{Proof of Theorem \ref{6}: 
the plan of campaign}\lab{campaign}

In this section we set out our framework for the proof of Theorem \ref{6}.
For any 
$(d,n) \in \mcal{E}$, where $\mcal{E}$ is given by (\ref{E}), 
let $F \in \Z[X_0,\ldots,X_n]$ be a non-singular form of degree
$d$.  Then the equation $F=0$ defines a non-singular hypersurface $X
\subset\bfP^{n}$ of degree $d$.
Recall that we have been following the convention that 
the implied constant in any estimate is absolute unless
explicitly indicated otherwise.  Since the pairs $(d,n)$ that occur
in the remainder of our work will always be restricted to lie in the set
$\mcal{E}$, we may henceforth assume that $d,n=O(1)$.

At the heart of our argument is an application of Lemma 
\ref{main1}.   According to this
result the points in which we are interested lie on 
$J=O_{\ve}(B^{n-1+\ve})$ linear subspaces
$M_1,\ldots, M_J \subset X$, all of which are defined over $\Q$.
Subspaces of dimension zero are clearly satisfactory for Theorem
\ref{6}, and so it remains to consider the case in which $M_j$ has
dimension $m \geq 1$.
In particular we know that $H(M_j)=O_{\ve}(B^{1+\ve})$ in these cases.
As an immediate  corollary of parts (i)--(iii) of Lemma \ref{fano}, it
follows that there are only $O(1)$ linear subspaces $M_j \subset X$
for which $m \geq(n-1)/2$.  Now Lemma \ref{triv} implies that any
single $m$-plane contributes $O(B^{m+1})$ to  $N_X(B)$, which is
satisfactory for Theorem \ref{6} since $m \leq n-2$.  Hence it remains
to handle the $m$-planes $M_j \subset X$ enumerated in Lemma
\ref{main1}, which have dimension 
\beq\lab{m-range}
1 \leq m < (n-1)/2.
\eeq
Since $n \leq 8$ for any $(d,n) \in \mcal{E}$, we plainly have
$m=1,2$ or $3$.

Let $1 \leq m < (n-1)/2$, and 
let $\Phi$ be an integral component of the Fano variety $F_m(X)$ of $m$-planes
contained in $X$.  When estimating the contribution from the
$m$-planes parametrised by $\Phi$, we shall always be able to 
assume that $X$ is a union of $m$-planes parametrised
by $\Phi$, in the sense that $X=D(\Phi)$ in the notation of (\ref{C}).
To see that this is permissible we suppose that the closed set $D(\Phi)$
is  a proper subvariety of $X$.  Then the degree of $D(\Phi)$ is
$O(1)$ by  part (i) of Lemma \ref{fano}, and so Lemma \ref{triv} 
shows that there is a 
contribution of $O(B^{n-1})$ to $N_X(B)$ from the points of
height at most $B$ that lie on $D(\Phi)$.  This is plainly
satisfactory for Theorem \ref{6}.  With this in mind we therefore
write 
$$
\tilde{F}_m(X) \subseteq F_m(X)
$$
for the union of integral components $\Phi \subseteq F_m(X)$ 
such that $D(\Phi)=X$.

When $n$ is odd, the possibility that $X$ may contain $(n-1)/2$-planes
is rather inconvenient, even though, as we have already seen, they 
make an acceptable contribution to $N_X(B)$.
We shall write $E$ for the union of such $(n-1)/2$-planes.
If $\Phi$ is an integral subvariety of
$F_m(X)$ we then write $\Phi^*$ for the open subset of $\Phi$ obtained by
removing any $m$-planes contained in $E$.  When
$n$ is even, so that $E$ is empty,  we merely take $\Phi^*=\Phi$.  

During the course of our work we shall have cause to estimate the
number $N_Y(H)$ of rational points of height at most $H$, on various
locally closed subsets $Y \subseteq X$, in the notation of
(\ref{flu}).  It will be convenient to combine here the information
that we gleaned in the previous section about large 
dimensional linear subspaces of $X$, together with
some of the estimates in \S \ref{prelim}.  Given any $k \in \N$, we define 
$$
\beta_k=\left\{
\begin{array}{ll}
0, & k \leq 3,\\
1/4, & k \geq 4.
\end{array}
\right.
$$
We then introduce functions $\sigma_n,\tau_n : \N \rightarrow
\Q$, given by
\beq\lab{sig}
\sigma_n(k)=\left\{
\begin{array}{ll}
k+1, & k\leq (n-1)/2,\\
k+\beta_k, & k>(n-1)/2,
\end{array}
\right. 
\eeq
and 
\beq\lab{tau}
\tau_n(k)=\left\{
\begin{array}{ll}
k+1, & k< (n-1)/2,\\
k+\beta_k, & k\geq (n-1)/2.
\end{array}
\right.
\eeq
When $n$ is even it is plain that $\sigma_n(k)=\tau_n(k)$ for each
$k \in \N$. 
The following result is a consequence of Lemmas
\ref{triv}, \ref{pila'} and \ref{fano}.

\begin{lem}\lab{glass}
Let $\ve>0$ and $H \geq 1$, and suppose that $Y \subseteq X$ is a subvariety of
dimension at most $k$, and degree $e$.  Then we have 
\begin{enumerate}
\item[(i)]
$
N_Y(H)\ll_{\ve,e} H^{\sigma_n(k)+\ve}.
$
\item[(ii)]
$
N_{Y\setminus E}(H)\ll_{\ve,e} H^{\tau_n(k)+\ve}.
$
\end{enumerate}
\end{lem}

Now suppose that $n=2m+1$ is odd, and let $k<m$.  Assume that
$\Phi\subseteq F_k(X)$ is integral, and that
$\Lambda\in F_{k-1}(X)$ is given.  Furthermore, let
$$
\Phi_{\Lambda}=\{\Gamma\in\Phi:\,\Lambda\subset \Gamma\}, 
\quad
\Phi_{\Lambda}^*=\{\Gamma\in\Phi^*:\,\Lambda\subset \Gamma\}.
$$
Then with this notation in mind, we proceed by establishing the following result.

\begin{lem}\lab{glass'}
Let $\ve>0$ and $H \geq 1$, and suppose that $D(\Phi_{\Lambda})$ has dimension
$\ell \leq m$ and degree $e$.  Then we have
\[
N_{D(\Phi_{\Lambda}^*)}(H)\ll_{\ve,e} H^{\tau_n(\ell)+\ve}.
\]
\end{lem}

\begin{proof}
Suppose that $\Pi\subseteq D(\Phi_{\Lambda})$ for
some $m$-plane $\Pi$, so that we must have $\ell=m$.  Since $D(\Phi_{\Lambda})$ is a cone with vertex
$\Lambda$ it follows that $\langle\Lambda,\Pi\rangle\subseteq
D(\Phi_{\Lambda})$, where $\langle\Lambda,\Pi\rangle$ is the linear space
spanned by $\Lambda$ and $\Pi$.  Since $\dim \Pi=\dim
D(\Phi_{\Lambda})=m$ this leads to a contradiction unless
$\Lambda\subseteq\Pi$. 

If $D(\Phi_{\Lambda})$ does not contain any $m$-planes, then the
statement of Lemma \ref{glass'} easily
follows from part (ii) of Lemma \ref{glass}.
Alternatively suppose that 
$D(\Phi_{\Lambda})$ contains an $m$-plane $\Pi$, say, so that in
particular $\ell=m$. 
Now let $x\in\Pi\subseteq D(\Phi_{\Lambda})$ be any point with
$x\not\in \Lambda$.  We claim that $x\not\in
D(\Phi_{\Lambda}^*)$.  But this follows from the observation that 
the $k$-plane 
$\langle x,\Lambda \rangle$ is contained in $\Pi$, so that 
$\langle x,\Lambda\rangle\not\in\Phi_{\Lambda}^*$, and the claim
follows. We have therefore shown that 
$$
N_{D(\Phi_{\Lambda}^*)}(H)\leq
N_{U}(H)+N_{\Lambda}(H),
$$ 
where $U \subseteq D(\Phi_\Lambda)$ denotes then open subset
formed by deleting all of the $m$-planes from $ D(\Phi_{\Lambda})$.
But then a simple application of part (ii) of Lemma~\ref{glass}
yields the required bound.  This completes the proof
of Lemma \ref{glass'}.
\end{proof}

Before embarking on the main thrust of the argument for Theorem \ref{6}
we take this opportunity to give an overview of the method.
Given an integral component $\Phi\subseteq \tilde{F}_m(X)$, the key idea
will be to estimate the contribution from the $m$-planes $\Lambda \in
\Phi$ defined over $\Q$, according to their smallest generator.  
Let $\Lambda \in \Grass(m,n)$ be any $m$-plane that is defined over
$\Q$. Then it is well known that $\Lambda$ contains 
linearly independent points $a_0,\ldots,a_m \in \bfP^n(\Q)$ such that
$$
H(a_0) \le \cdots \le H(a_m), \quad H(\Lambda) \ll \prod_{i=0}^m
H(a_i) \ll H(\Lambda).
$$
We shall call a rational point on $\Lambda$ a 
``smallest generator'' for $\Lambda$ if
it has minimal height.  Any smallest generator 
for $\Lambda$ obviously has height
$O(H(\Lambda)^{1/(m+1)})$.

We now turn to the pivotal role that Lemma \ref{strat} plays in this work.
In fact it helps us to estimate the contribution
from the $m$-planes in $\Phi \subseteq \tilde{F}_m(X)$ 
according to their smallest generator, $a_0$, say, in essentially
two different ways.  The first approach involves counting the
total number of rational points of height at most $B$ that lie on the
cone $C_{a_0}(\Phi)=D(\Phi_{a_0})$. (Recall that this is the cone 
swept out by the $m$-planes contained in $\Phi$ and passing through $a_0$.)  
This will be referred to as ``counting by cones'', or the
``\cc-method'' for short.  In the second approach one counts the
total number of $m$-planes of height $O_\ve(B^{1+\ve})$ that pass
through $a_0$, before then estimating the number of rational points of
height at most $B$ on each such $m$-plane.  
This will be referred to as ``counting by linear spaces'', or the
``\cl-method'' for short.
In both approaches one then obtains a final estimate for the
contribution to $N_X(B)$, by summing over all of the points 
$a_0 \in X(\Q)$ of low height.
In doing so we shall need to apply Lemma \ref{strat} with the 
choice $Y=X$, in order to be
able to control the dimension of $\Phi_{a_0}$ as $a_0$ varies in $X$.

We take a moment to analyse the \cl-method in more detail.  It will
become apparent that our implementation of this approach has a distinctly
subtler nature than that of the \cc-method.  
We shall focus upon the case $m=2$ here, the same principle applying
in simpler form for the case $m=1$.  Let
\beq\lab{range}
A_2\geq A_1 \geq A_0 \geq 1,
\eeq
and let $\Phi \subseteq
\tilde{F}_2(X)$ be an integral component as above.
The method begins by fixing a point $a_0
\in X(\Q)$, which has height $A_0/2<H(a_0) \leq A_0$.
One then considers the possible rational points $a_1$ in the
cone $C_{a_0}(\Phi)$, which have height $A_1/2<H(a_1) \leq A_1$.
This then fixes the smallest two generators of the planes we wish to
estimate the contribution from.  Finally, for fixed $a_0,a_1$ one
forms the variety 
$$
\Phi_{a_0,a_1}=\{P \in \Phi_{a_0}: a_1 \in P\}
$$ 
of planes in $\Phi$ which contain the line $\langle{a_0,a_1}\rangle$,
and then estimates the number of $a_2\in D(\Phi_{a_0,a_1})$ such that
$A_2/2<H(a_2) \leq A_2$.   Each such value of $a_2$ determines a plane
$P=\langle{a_0,a_1,a_2}\rangle\in \Phi$, which it may be assumed has height
$$
A_0A_1A_2\ll H(P) \ll A_0A_1A_2.
$$
One then employs Lemma \ref{triv-lin}
to estimate the number of rational points of height at
most $B$ contained in $P$.
We need to control the dimension of 
$\Phi_{a_0}$ as $a_0$ varies in $X$.  Similarly, for fixed $a_0$, we
need to control the 
dimension of $\Phi_{a_0,a_1}$
as $a_1$ varies in $C_{a_0}(\Phi)$.  In the first case we shall apply
Lemma \ref{strat} with $Y=X$ as indicated above,
and in the second we shall apply the same lemma, but with $Y=D(\Psi)$,
where $\Psi$ is an integral component of $\Phi_{a_0}$.
Ultimately, since one is only interested in planes of height
$O_{\ve}(B^{1+\ve})$, one sums over dyadic intervals for
$A_0,A_1,A_2$ such that (\ref{range}) holds and $A_0A_1A_2
\ll_{\ve} B^{1+\ve}$.
When $m=1$ the \cl-method is plainly simpler since we need 
only apply Lemma \ref{strat} once.

The astute reader will notice that there is a degree of waste in the
above description of the \cl-method.  
Thus for fixed $a_0 \in  X(\Q)$, whereas we are only interested in the
number of lines of height $A_0A_1$ that pass through $a_0$, we are in
fact counting the total number of possible $a_1$ that serve as
generators for the line $\langle{a_0,a_1}\rangle$ of height $A_0A_1$.
In fact Lemma \ref{triv-lin} implies that any such line contains
$\gg A_1/A_0$ rational points $a_1$ such that $H(a_1) \leq A_1$.  
Hence it follows that the total number of rational lines of height
$A_0A_1$ contained in $\Phi_{a_0}$ is actually 
\beq\lab{over-1}
\ll \frac{A_0}{A_1} \#\{a_1 \in C_{a_0}(\Phi)(\Q): H(a_1) \leq A_1\},
\eeq
in which the first factor allows for the ``over-counting'' inherent in our
method.  A similar phenomenon occurs when $a_0,a_1$ are 
fixed and one is counting the number of planes of height
$A_0A_1A_2$ that pass through the line $\langle{a_0,a_1}\rangle$.
Thus an application of Lemma \ref{triv-lin} reveals that 
each such plane contains $\gg A_2^2/A_0A_1\gg A_2/A_0$ suitable points
$a_2$, whence the total number of planes of height
$A_0A_1A_2$ contained in $\Phi_{a_0,a_1}$ is
\beq\lab{over-2}
\ll \frac{A_0}{A_2} \#\{a_2 \in D(\Phi_{a_0,a_1})(\Q): H(a_2) \leq A_2\}.
\eeq
In summary the \cl-method has two main ingredients: 
a stratification argument involving Lemma~\ref{strat}, 
and a means of rectifying the
over-counting that our implementation of the \cl-method engenders.

We now have all of the tools with which to complete the proof of
Theorem~\ref{6}.
At this stage it is convenient to make a certain 
hypothesis concerning the quantity $N_X(B)$.

\begin{hyp}[Projective hypersurface hypothesis]
Let $\ve>0$ and let $d,n \geq 2$ be integers.
Then there exists $\theta_{d,n}\geq 0$ such that
$$
N_X(B)=O_{d,\ve,n}(B^{n-1+\theta_{d,n}+\ve}),
$$
for any non-singular hypersurface $X \subset \bfP^n$ of degree $d$, that is defined
over $\Q$.
\end{hyp}

We shall henceforth write $\ph[\theta_{d,n}]$ to denote the 
projective hypersurface
hypothesis holding with exponent $\theta_{d,n}$.
Thus Conjecture \ref{hb-ns} is the
statement that $\ph[0]$ holds, and it follows from Lemma \ref{pila'}
that $\ph[1/4]$ holds.  For our purposes it will actually suffice to note that
$\ph[1]$ holds, by Lemma \ref{triv}.

For any integer $m$ in the range (\ref{m-range}), we let $X_m \subseteq
X$ be the finite union of $m$-planes $M_j \subset X$ that are enumerated in
Lemma \ref{main1}, and that are parametrised by $\tilde{F}_m(X)$.  
In the notation of (\ref{flu}), we shall write $N_{X_m}(B)$ for the
overall contribution to $N_X(B)$ from the points contained in $X_m$.
Our main task is to prove the following result.

\begin{pro}\lab{Xm}
Let $(d,n) \in \mcal{E}$ and let $1
\leq m < (n-1)/2$.  Then we have 
$$
N_{X_m}(B) \ll_{\ve} B^{n-1+\theta_{d,n}/2+\ve},
$$
provided that $\ph[\theta_{d,n}]$ holds. 
\end{pro}

Before establishing Proposition \ref{Xm}, we first indicate how it can be
used to prove Theorem \ref{6}.
Suppose that $\ph[\theta]$ holds for $\theta=\theta_{d,n}\geq 0$. 
Then it easily follows
from Proposition \ref{Xm} and our work in the previous section that 
$$
N_X(B) \ll_{\ve} B^{n-1+\theta/2+\ve}.
$$
But then we see that $\ph[\theta/2]$ holds. 
This allows us to deduce the sharper bound 
$$
N_X(B) \ll_{\ve} B^{n-1+\theta/4+\ve},
$$
whence in fact $\ph[\theta/4]$ holds.   
By continuing to iterate this procedure sufficiently many
times we may clearly conclude that $\ph[\ve]$ holds for any given
$\ve>0$. This completes  the proof of Theorem \ref{6},
upon revising the choice of $\ve$.

For the remainder of this paper we shall assume that
$\ph[\theta]$ holds for some $\theta\geq 0$, and that this
value of $\theta$ is identical for each $(d,n)\in \mcal{E}$.
Furthermore, we shall follow common practice and allow the small
positive constant $\ve$ to take different values at different parts of
the argument.

\section{Proof of Theorem \ref{6}: lines}
\lab{line}

We begin the proof of Proposition \ref{Xm} by handling the contribution from
the lines $M_j \subset X_1$, for which we shall employ the \cl-method
that was introduced in \S \ref{campaign}.
Let $\Phi\subseteq\tilde{F}_1(X)$ be any integral
component.  On applying Lemma \ref{strat} with $Y=X$ we obtain a 
stratification of subvarieties 
$X=Z_0(\Phi) \supset Z_1(\Phi)$ 
such that $\deg Z_1(\Phi)=O(1)$ and 
\beq\lab{s-i}
\dim Z_1(\Phi) \leq n-3, 
\eeq 
and 
\beq\lab{s-ii}
\dim \Phi_y = \dim \Phi-n+2,
\eeq
for any $y \in X\setminus Z_{1}(\Phi)$.

Considering the component $\Phi \subseteq \tilde{F}_1(X)$ as being
fixed, it will be convenient to write $Z_1=Z_1(\Phi)$.
Our plan will be to sort the lines in $\Phi$ according to whether 
their smallest generator lies in $X\setminus Z_1$, or in $Z_1$. 
Thus we shall write $M_0(B)$ for the overall contribution to
$N_{X_1}(B)$ from the  lines contained in $X_1$ that are parametrised by
$\Phi$ and have smallest generator $a_0 \in X \setminus Z_1$, 
and we shall write $M_1(B)$ for the corresponding contribution from
the lines with smallest generator $a_0 \in Z_1$. 
Now there are $O(1)$ possible integral components of $\tilde{F}_1(X)$.
In order to establish Proposition \ref{Xm} in the case $m=1$ it will
therefore suffice to show that
$$
M_i(B)=O_\ve (B^{n-1+\theta/2+\ve}),
$$
for $i=0,1$.  Recall that any line contained in $X_1$
has height $O_{\ve} (B^{1+\ve}).$ For any real numbers $A_0,A_1$ such
that
\beq\lab{01}
A_1\geq A_0 \geq 1, \quad A_0A_1 \ll_{\ve}B^{1+\ve},
\eeq
let  $M_i(B;A_0,A_1)$ denote the contribution to
$M_i(B)$ from the lines of height $A_0A_1$ that pass
through rational points $a_0$ and $a_1$, such that 
$a_0$ is a smallest generator for the line and
$$
A_0/2< H(a_0)  \leq A_0, \quad A_1/2< H(a_1) \leq A_1.
$$
On summing over $O_{\ve}(B^{\ve})$ dyadic intervals for $A_0,A_1$,
it will therefore suffice to show that
\beq\lab{goal-i}
M_i(B;A_0,A_1)=O_\ve (B^{n-1+\theta/2+\ve}),
\eeq
for $i=0,1$, and each choice of $A_0,A_1$ such that (\ref{01}) holds.

We begin by establishing (\ref{goal-i}) in the case $i=0$.
For any $a_0 \in X\setminus Z_1$, recall the definitions 
of the cones $\Phi_{a_0}$ and
$C_{a_0}(\Phi)=D(\Phi_{a_0})$, as given by (\ref{C}) and
(\ref{fibred}).
Since $(d,n) \in \mcal{E}$, it follows from part (iv) of Lemma
\ref{fano} that $\Phi$ has dimension at most $2n-3-d$. Hence, we
deduce from (\ref{s-ii}) that 
$$
\dim C_{a_0}(\Phi) \leq n-d.
$$
Moreover we have already seen that $\deg C_{a_0}(\Phi)=O(1)$.
We claim that
$$
\#\{a_1 \in C_{a_0}(\Phi^*)(\Q): H(a_1) \leq A_1\} =O_{\ve}(A_1^{n-3+\ve}),
$$
where $C_{a_0}(\Phi^*)$ is the union of those lines in $\Phi^*$ that
pass through $a_0$.
This follows from Lemma \ref{triv} when $d=4$.
If $d=3$ and $n \geq 7$ we know from Lemma \ref{deg-cone} that
$C_{a_0}(\Phi)=C_{a_0}^1$ is
integral and has degree $6$, so that the result follows from 
Lemma \ref{pila'}.  Finally for the case $(d,n)=(3,6)$, we apply
part (ii) of Lemma \ref{glass}, and for the case $(d,n)=(3,5)$ we apply Lemma 
\ref{glass'}.  

In order to estimate the total number of rational lines of height
$A_0A_1$ that are parametrised by $\Phi_{a_0}^*$, we employ the
over-counting argument used in (\ref{over-1}) to deduce that there are
$$
\ll \frac{A_0}{A_1} \#\{a_1 \in C_{a_0}(\Phi^*)(\Q): H(a_1) \leq A_1\}
\ll_\ve \frac{A_0}{A_1}A_1^{n-3+\ve} = A_0A_1^{n-4+\ve}
$$
such lines.  Moreover Lemma \ref{triv-lin} implies that 
any line of height $A_0A_1$ contains $O_{\ve}(B^{2+\ve}/(A_0A_1))$ points of
height at most $B$, since $A_0A_1 \ll_{\ve}B^{1+\ve}$ by (\ref{01}).
 Putting all of this together we therefore obtain 
\begin{align*}
M_0(B;A_0,A_1)&\ll_{\ve} B^{n-1} +\sum_{\colt{a_0 \in
    (X\setminus Z_1)(\Q)}{H(a_0)  \leq A_0}} 
A_0A_1^{n-4+\ve}.\frac{B^{2+\ve}}{A_0A_1}\\
&\ll_\ve B^{n-1}+ A_0^{n-1+\theta}
A_1^{n-5}B^{2+\ve}.
\end{align*}
Thus  (\ref{01}) yields
$$
M_0(B;A_0,A_1)
\ll_{\ve}  B^{n-1}+ A_0^{4+\theta}B^{n-3+\ve}
\ll_{\ve} B^{n-1+\theta/2+\ve},
$$
which thereby establishes (\ref{goal-i}) in the case $i=0$.

We now turn to the proof of (\ref{goal-i}) for $i=1$.
Let $A_0,A_1$ be real numbers such that (\ref{01}) holds.
For any $a_0\in Z_1$, it follows from (\ref{condition}) that 
$\dim C_{a_0}(\Phi) \le n-2$, and from  (\ref{s-i}) 
that $Z_1$ has dimension at most $n-3$.
Recall the definitions (\ref{sig}) and (\ref{tau}) of $\sigma_n$ and
$\tau_n$, respectively.  Then part (i) of Lemma \ref{glass} implies that
$Z_1$ contains $O_{\ve}(A_0^{\sigma_n(n-3)+\ve})$ rational points of height
$A_0$, and that $C_{a_0}(\Phi)$ 
contains $O_{\ve}(A_1^{\sigma_n(n-2)+\ve})$ rational points of height
$A_1$.  We therefore obtain the estimate
\begin{align*}
M_1(B;A_0,A_1)
&\ll_{\ve}  B^{n-1}+ 
A_0^{\sigma_n(n-3)}.\frac{A_0}{A_1}.A_1^{\sigma_n(n-2)}
\frac{B^{2+\ve}}{A_0A_1}\\
&\ll_{\ve}  B^{n-1}+ 
A_0^{\sigma_n(n-3)} A_1^{\sigma_n(n-2)-2}B^{2+\ve},
\end{align*}
for $i \geq 1$.
Suppose first that $n=5$, so that 
$\sigma_n(n-3)=\sigma_n(n-2)=3$.  Then it follows from
(\ref{01}) that
$$
M_1(B;A_0,A_1)
\ll_{\ve}  B^{4}+ 
A_0^{3} A_1 B^{2+\ve}
\ll_{\ve}  B^{4+\ve},
$$
which is satisfactory for (\ref{goal-i}).  
If however $n\geq 6$, then 
$$
\sigma_n(n-3)\leq n-11/4, \quad 
\sigma_n(n-2)= n-7/4,
$$
and we see that for sufficiently small $\ve>0$ we have
$$
M_1(B;A_0,A_1)
\ll_{\ve}  B^{n-1}+ 
A_0^{n-11/4} A_1^{n-15/4}B^{2+\ve} \ll B^{n-1}.
$$
This completes the proof of
(\ref{goal-i}), and so the proof of
Proposition \ref{Xm} for  $m=1$.

\section{Proof of Theorem \ref{6}: planes}
\lab{plane}

Next we consider the case $m=2$ of planes $M_j\subset X_2$, 
for $(d,n) \in \mcal{E}$.  In
particular we may henceforth assume that $n\geq 6$, since $m<(n-1)/2$
by (\ref{m-range}).  We begin by dispatching the case $n=6$, for which
we may assume that $d=4$.  Indeed 
we have already observed in the context of
Conjecture \ref{hope} that a non-singular cubic 
hypersurface $X \subset \bfP^6$ is not a
union of planes, so that in particular $\tilde{F}_2(X)$ is empty.
Assuming therefore that $(d,n)=(4,6)$, and that $\tilde{F}_2(X)$ is
non-empty, we proceed by employing a
simple version of the \cl-method that was outlined in \S \ref{campaign}.
The planes $P$ in which we are interested have height
$H(P)=O_\ve(B^{1+\ve})$, so any smallest generator for $P$ must
have height $O_\ve(B^{1/3+\ve})$.
Now it follows from Lemma \ref{fano'} that $\tilde{F}_2(X)$ has
dimension at most $3$.
If $\Phi\subseteq\tilde{F}_2(X)$ is any integral component 
then Lemma \ref{strat} yields a 
stratification $X=Z_0(\Phi) \supset Z_1(\Phi),$
such that $\deg Z_1(\Phi)=O(1)$, $\dim Z_1(\Phi) \leq 3$, and
\beq\lab{r-ii}
\dim \Phi_y = \dim \Phi-5+2 \leq 0
\eeq
for any $y \in X \setminus Z_1(\Phi)$.
We can and will assume that the component $\Phi$ is fixed,
so that we may write $Z_1=Z_1(\Phi)$ for convenience.

We shall write $M_0(B)$ for the overall contribution to
$N_{X_2}(B)$, arising from the planes whose smallest generator lies in
$X \setminus Z_1$. In this setting it follows from (\ref{r-ii}) 
that each point $a_0\in X\setminus Z_1$ gives rise to only $O(1)$ planes
in $\Phi_{a_0}$, and each such plane contributes $O(B^3)$ by
Lemma \ref{triv}.  Thus a further application of Lemma \ref{triv} yields
$$
M_0(B)\ll \sum_{\colt{a_0 \in
    X(\Q)}{H(a_0)\ll_{\ve}B^{1/3+\ve}}} B^3 \ll_\ve 
(B^{1/3+\ve})^6 B^3 \ll_\ve B^{5+\ve},
$$
which is satisfactory for Proposition \ref{Xm}. 
Next we tackle the contribution $M_1(B)$ from the planes that 
have smallest generator $a_0 \in
Z_1$.  Lemma \ref{triv-lin} implies that
$$
M_1(B)\ll_\ve \sum_{\colt{a_0 \in
    Z_1(\Q)}{H(a_0)\ll_{\ve}B^{1/3+\ve}}}
\sum_{R}
\sum_{\colt{P \in \Phi_{a_0}}{R/2<H(P)\leq R}}
\frac{B^{3+\ve}}{R},
$$
where the summation over $R$ is over $O_\ve(B^\ve)$ dyadic intervals
for $R \ll_\ve B^{1+\ve}$.
We have seen that $Z_1$ has dimension at most $3$,  and
we claim that $\Phi_{a_0}$ has dimension at most $2$ for any $a_0 \in
Z_1$. To see this we note that if $\Phi_{a_0}\subseteq \Phi$
had dimension $3$ for some $a_0 \in X$, then it would follow that
$\Phi_{a_0}=\Phi$, since $\Phi$ is integral and has 
dimension at most $3$.  
However $D(\Phi_{a_0})$ is a cone with vertex $a_0$, so that
$D(\Phi_{a_0})$ cannot be a non-singular hypersurface of 
degree $2$ or more.  Since $D(\Phi)=X$ is non-singular we obtain a 
contradiction, which establishes the claim. 
In view of part (ii) of Lemma \ref{fano} we see that $\Phi_{a_0}$
does not contain any linear spaces of dimension $2$.  
Indeed, by \cite[p.123]{harris}, even a
line in $F_2(X)$ would produce a collection of planes which spanned a
$3$-plane, contradicting Lemma \ref{fano}, part (ii).  We also see that
$Z_1$ does not contain any $3$-planes. Hence Lemma
\ref{pila'} implies that 
$$
M_1(B)\ll_\ve \sum_{\colt{a_0 \in
    Z_1(\Q)}{H(a_0)\ll_{\ve}B^{1/3+\ve}}}
\sum_{R \ll_\ve B^{1+\ve}}
R^{1+\ve}B^{3+\ve}\ll_\ve 
(B^{1/3+\ve})^{3+\ve}B^{4+\ve} \ll_\ve B^{5+\ve}.
$$
This too is satisfactory for Proposition \ref{Xm}, and so completes
the proof of this result in the case $m=2$ and $n=6$.

For the rest of this section we shall assume that $n \geq 7$ so that
together Lemmas~\ref{fano} and \ref{fano'} imply that
\beq\lab{f2}
\dim \tilde{F}_2(X) \leq 3n-16.
\eeq
Let $\ma{A}=(A_0,A_1,A_2)$, with 
\beq\lab{range-Ai}
A_2\geq A_1 \geq A_0 \geq 1, \quad A_0A_1A_2 \ll_\ve
B^{1+\ve}.
\eeq  
Then our objective is to estimate the contribution
$N_{X_2}(B;\ma{A})$, say,  from the planes $P \in \tilde{F}_2(X)$ with height
of order $A_0A_1A_2$, that are generated by linearly independent 
points $a_0, a_1, a_2 \in X(\Q)$ such that 
\beq\lab{height-stu}
A_0/2< H(a_0)  \leq A_0, \quad A_1/2< H(a_1) \leq A_1, \quad
A_2/2< H(a_2) \leq A_2.
\eeq
In doing so it will plainly suffice to estimate 
the contribution from those planes that belong to $\Phi$, for a
fixed integral component $\Phi \subseteq \tilde{F}_2(X)$.
Our plan will be to apply the \cl-method, as described in \S \ref{campaign}.

Suppose that $a_0 \in X(\Q)$ is a fixed point such that $A_0/2< H(a_0)
\leq A_0$, and let $\Psi$ be an integral component of $\Phi_{a_0}$.
Then we are interested in the planes that are
parametrised by $\Psi$, the union of which form a cone $D(\Psi)$.
We claim that it suffices to assume that
\beq\lab{woodstock}
\dim D(\Psi) \geq n-3.
\eeq
To see that this is
permissible, we suppose for the moment that $D(\Psi)$ has dimension at
most $n-4$, and deduce from part (ii) of Lemma  \ref{glass} 
that there are $O_\ve(B^{\tau_n(n-4)+\ve})$
rational points of height at most $B$ contained in $D(\Psi)\setminus E$.
Since $n\ge 7$ we have $\tau_n(n-4)=n-4+1/4$.
Moreover there are $O(A_0^n)$ possible choices for $a_0$ by Lemma
\ref{triv}. On applying (\ref{range-Ai}) we therefore obtain the
overall contribution
$$
\ll_\ve (B^{1/3+\ve})^nB^{n-4+1/4+\ve} \ll_\ve B^{4n/3-4+1/4+\ve}
$$
to $N_{X_2}(B,\ma{A})$ from this scenario.  Once summed over dyadic
intervals for the $A_0,A_1,A_2$ this is plainly satisfactory for
Proposition \ref{Xm}, since $n \leq 8$. 
In fact, strictly speaking,
summation over dyadic values of $A_1,A_2$ is unnecessary, since the above bound deals
with all planes through $a_0$. It therefore suffices to assume
that (\ref{woodstock}) holds for each $a_0 \in X(\Q)$ and each integral
component $\Psi$ of $\Phi_{a_0}$.

We now apply Lemma \ref{strat} with $Y=D(\Psi)$. This produces a
stratification of subvarieties 
$$
D(\Psi)=W_0(\Psi) \supseteq W_1(\Psi) \supseteq W_2(\Psi) \supseteq \cdots,
$$
such that $\deg W_i(\Psi)=O(1)$, with 
\beq\lab{t-i}
\dim W_i(\Psi) \leq \dim D(\Psi)-1-i,\quad (i\ge 1),
\eeq
and 
\beq\lab{t-ii}
\dim \Psi_y \leq \dim \Phi_{a_0}-\dim D(\Psi)+2+i,\quad (i\ge 0),
\eeq
for any $y \in W_i(\Psi) \setminus W_{i+1}(\Psi)$. 
We proceed to estimate the contribution from those planes 
$P \in \Phi_{a_0}$, with height of order
$A_0A_1A_2$, which are generated by linearly independent rational points
$a_0, a_1, a_2$ such that (\ref{height-stu}) holds.  We shall do this 
according to the value of $i \geq 0 $ for which
$$
a_1 \in W_i(\Psi) \setminus W_{i+1}(\Psi).
$$

For each fixed value of $a_1$, we have $a_2 \in
C_{a_1}(\Psi)=D(\Psi_{a_1})$, which we suppose has dimension
$\alpha_2=\alpha_2(i; a_0,a_1,\Psi)$. In particular we must have $\al_2
\geq 2$ whenever $\Psi_{a_1}$ is non-empty.  Taken together with Lemma
\ref{glass'}, part (ii) of Lemma \ref{glass} now shows that there are
$O_{\ve}(A_2^{\tau_n(\alpha_2)+\ve})$ relevant points $a_2$,
where $\tau_n$ is given by (\ref{tau}).  By the over-counting argument used in
(\ref{over-2}) this produces
$O_{\ve}(A_0A_2^{\tau_n(\alpha_2)-1+\ve})$ 
planes of height $A_0A_1A_2$ through the line $\langle{a_0,a_1}\rangle$.
In the special case $\alpha_2=2$ we may improve this, since
there are $O(1)$ planes in $\Psi_{a_1}$.  We therefore write
\beq\lab{mu}
\mu_n(k)=\left\{
\begin{array}{ll}
\tau_n(k), & k\geq 3,\\
1, & k\leq 2,
\end{array}
\right.
\eeq
and deduce that there are
$\ll_\ve A_0 A_2^{\mu_n(\alpha_2)-1+\ve}$
available planes.

Now write $\alpha_1=\alpha_1(i; a_0,\Psi)$ for the dimension of
$W_i(\Psi)$, so that  Lemma \ref{glass}, part (i),  shows there to be 
$O_{\ve}(A_1^{\sigma_n(\alpha_1)+\ve})$
available points $a_1 \in W_i(\Psi) \setminus W_{i+1}(\Psi)$, where
$\sigma_n$ is given by (\ref{sig}).  This time the over-counting 
argument used in
(\ref{over-1}) shows that there are 
$$
\ll_\ve \frac{A_0}{A_1}. A_1^{\sigma_n(\alpha_1)+\ve}
= A_0 A_1^{\sigma_n(\alpha_1)-1+\ve}
$$
available lines corresponding to points 
$a_1 \in W_i(\Psi) \setminus W_{i+1}(\Psi)$.
In conclusion we have therefore shown that for fixed $a_0 \in X$, and any fixed
integral component $\Psi \subseteq \Phi_{a_0}$, the overall number of
planes in $\Phi_{a_0}$ that have height of order $A_0A_1A_2$, 
and whose smallest generators are of order $A_0,A_1,A_2$,
respectively, is
\beq\lab{tar}
\ll_\ve A_0^2 A_1^{\sigma_n(\alpha_1)-1}A_2^{\mu_n(\alpha_2)-1+\ve},
\eeq
for certain integers $\al_1,\al_2 \geq 0$.

Before going on to consider the corresponding number of points $a_0$,
we record some useful inequalities concerning the quantities $\al_1$
and $\al_2$ introduced above.  Suppose first that $i \geq 1$ in the
stratification.  Then it follows from (\ref{t-i}) and (\ref{t-ii})
that
\beq\lab{i>0}
\al_1+\al_2 \leq  \dim \Phi_{a_0}+3, \quad (i\geq 1).
\eeq
Alternatively, if $i=0$ then we have $D(\Psi_{a_1}) \subseteq
D(\Psi)=W_0(\Psi)$, leading to the inequality $\al_2 \leq \al_1$.  Thus it 
follows from (\ref{woodstock}) and (\ref{t-ii}) that 
\beq\lab{i=0}
2\leq \al_2,n-3\leq \al_1, \quad \al_1+\al_2 \leq \dim \Phi_{a_0}+4,
\quad (i=0),
\eeq
in this case. Here, as throughout the remainder of our work, the first
set of inequalities is always taken to mean
$2 \leq \min\{\al_2,n-3\} \leq \max\{\al_2,n-3\} \leq \al_1$.
Finally we record the trivial inequalities 
\beq\lab{al-triv}
\al_1 \leq n-2, \quad \al_2 \leq n-3, \quad (i\ge 0), 
\eeq
that follow from (\ref{condition}) and Lemma \ref{condition'}.
It is important to note that the inequalities 
(\ref{i>0})--(\ref{al-triv}) are valid for any
choice of $a_0 \in X(\Q)$.
Moreover, there are clearly $O(1)$ possible integral components $\Psi
\subseteq \Phi_{a_0}$. Hence we deduce that for fixed $a_0 \in
X$ the number of planes in $\Phi_{a_0}$ that have height
of order $A_0A_1A_2$ is given by (\ref{tar}), where $\al_1,\al_2 \geq 0$ 
satisfy (\ref{i>0}) or
(\ref{i=0}), together with (\ref{al-triv}).

We now turn to the problem of estimating the number of possible smallest 
generators $a_0$.  Here we shall
combine (\ref{f2}) with an application of Lemma \ref{strat} in the
case $Y=X$.  Thus there exists a 
stratification of subvarieties 
$$
X=Z_0(\Phi) \supseteq Z_1(\Phi) \supseteq Z_2(\Phi) \supseteq \cdots,
$$
such that $\deg Z_j(\Phi)=O(1)$, with
\beq\lab{s-j}
\dim Z_j(\Phi) \leq n-2-j,\quad(j\ge 1),
\eeq 
and 
\beq\lab{s-jj}
\dim \Phi_y \leq 2n-13+j, \quad (j\ge 0),
\eeq
for any $y \in Z_j(\Phi) \setminus Z_{j+1}(\Phi)$.
We write $Z_j=Z_j(\Phi)$ for $j \geq 0$, for convenience.
In keeping with the above,  our plan is to estimate the overall
contribution from the planes  
$P \in \Phi$ with height of order $A_0A_1A_2$, that are generated by 
points $a_0, a_1, a_2$ such that
(\ref{height-stu}) holds.  We shall classify such planes 
according to the value of $j \geq 0 $ for
which $a_0 \in Z_j\setminus Z_{j+1}.$
Write $\alpha_0=\alpha_0(j)$ for the dimension of $Z_j$, and let
\beq\lab{nu}
\nu_{n,\theta}(k)=\left\{
\begin{array}{ll}
\sigma_n(k), & k\leq n-2,\\
k+\theta, & k=n-1,
\end{array}
\right.
\eeq
for $k \in\N$.  Then it follows from part (i) of Lemma \ref{glass} and the fact
that $\ph[\theta]$ holds, that the total number of available points $a_0$ is
$$
\ll_\ve A_0^{\nu_{n,\theta}(\alpha_0)+\ve}.
$$
On combining this with (\ref{tar}) we therefore conclude that the
total number of planes under consideration is
\beq\lab{a0}
\ll_\ve A_0^{\nu_{n,\theta}(\alpha_0)+2}A_1^{\sigma_n(\alpha_1)-1}
A_2^{\mu_n(\alpha_2)-1+\ve},
\eeq
for certain integers $\al_0,\al_1,\al_2\geq 0$.

We now collect together some of the inequalities satisfied by
$\al_0,\al_1,\al_2$ as we range over values of $i,j \geq 0$ in our
double stratification.
Suppose first that $j \geq 1$.  Then it follows from (\ref{s-j}) that
$\al_0 \leq n-2-j$, whence we may combine (\ref{i>0}), (\ref{i=0}) and
(\ref{s-jj}) to deduce that 
\beq\lab{final-1}
\al_0+\al_1+\al_2 \leq 3n-12, \quad (i,j\ge 1),
\eeq
and 
\beq\lab{final-2}
2\leq \al_2,n-3 \leq \al_1, \quad \al_0+\al_1+\al_2 \leq 3n-11, \quad (i=0,\, j\ge 1).
\eeq
Similarly, since $Z_0=X$, we see that if $j=0$ then 
(\ref{i>0}), (\ref{i=0}) and (\ref{s-jj}) combine to give
\beq\lab{final-3}
\al_0=n-1, \quad \al_0+\al_1+\al_2 \leq 3n-11, \quad (i\ge 1,\,j=0), 
\eeq
and 
\beq\lab{final-4}
\al_0=n-1, \quad 2\leq \al_2,n-3\leq \al_1, 
\quad \al_0+\al_1+\al_2 \leq 3n-10, \quad (i=j=0). 
\eeq

We are now ready to complete our treatment of the planes in
Proposition~\ref{Xm}.
By Lemmas \ref{triv} and \ref{triv-lin} it follows that any plane of height
$A_0A_1A_2$ contains
$$
\ll B^2 + \frac{B^3}{A_0A_1A_2} \ll_\ve \frac{B^{3+\ve}}{A_0A_1A_2}
$$
rational points of height at most $B$.  On combining this with
(\ref{a0}), we therefore conclude the proof of the following result.

\begin{lem}\lab{oriel}
Suppose that $A_0,A_1,A_2$ satisfy (\ref{range-Ai}). 
Then there exists a triple $\mal=(\al_0,\al_1,\al_2)$ of non-negative
integers satisfying one of the conditions
(\ref{final-1})--(\ref{final-4}), together with (\ref{al-triv}), such that
$$
N_{X_2}(B;\ma{A}) \ll_\ve 
A_0^{\nu_{n,\theta}(\alpha_0)+1}A_1^{\sigma_n(\alpha_1)-2}
A_2^{\mu_n(\alpha_2)-2}B^{3+\ve}.
$$
Here $\nu_{n,\theta}, \sigma_n, \mu_n$ are given by (\ref{nu}),
(\ref{sig}) and (\ref{mu}), respectively.
\end{lem}

It remains to apply the linear programming result Lemma \ref{prog} in
order to show that the bound in Lemma
\ref{oriel} is satisfactory for Proposition \ref{Xm}. 
In view of the inequalities (\ref{range-Ai}) satisfied by
$A_0,A_1,A_2$, we deduce from Lemmas \ref{prog} and \ref{oriel} that
$$
N_{X_2}(B;\ma{A}) \ll_\ve \Big(M_1(B)+M_2(B)\Big)B^{3+\ve},
$$
where
$$
M_1(B)=B^{\mu_n(\alpha_2)-2}+
B^{(\sigma_n(\alpha_1)+\mu_n(\alpha_2)-4)/2},
$$
and 
$$
M_2(B)=B^{(\nu_{n,\theta}(\alpha_0)+\sigma_n(\alpha_1)+\mu_n(\alpha_2)-3)/3}.
$$
On summing over dyadic intervals for the $A_0,A_1,A_2$, we see that it
suffices to show that 
$$
M_i(B)=O(B^{n-4+\theta/3}),
$$
for $i=1,2$, in order to complete the proof of Proposition \ref{Xm} in
the case $m=2$.
We first establish the case $i=1$ of this estimate.  
On recalling the definitions
(\ref{sig}) and (\ref{mu}) of $\sigma_n$ and $\mu_n$, respectively, we
may clearly apply (\ref{al-triv}) to deduce that
$$
M_1(B)\leq B^{\mu_n(n-3)-2}+
B^{(\sigma_n(n-2)+\mu_n(n-3)-4)/2}
\leq B^{n-19/4}+
B^{n-17/4}.
$$
This is plainly satisfactory.

Turning to the estimate for $M_2(B)$, we have four different cases to
consider.  In each one our task is to show that 
\beq\lab{En}
E_n(\mal)=\frac{\nu_{n,\theta}(\alpha_0)+\sigma_n(\alpha_1)
+\mu_n(\alpha_2)}{3}-1\leq n-4+\theta/3.
\eeq
Suppose firstly that the triple $\mal$
satisfies (\ref{final-1}).  Then it is trivial to see that
$E_n(\mal) \leq (3n-9)/3-1=n-4,$
which is satisfactory for (\ref{En}).
Next suppose that $\mal$ satisfies
(\ref{final-3}).  Then $\al_0=n-1$ and $\al_1+\al_2\leq 2n-10$, from
which it follows that
$E_n(\mal)\leq (3n-9+\theta)/3-1=n-4+\theta/3$.
This too is satisfactory for (\ref{En}).
In order to handle the remaining two cases in which $\mal$
satisfies (\ref{final-2}) or (\ref{final-4}), 
it plainly suffices to assume that $\al_0+\al_1+\al_2 \geq 3n-11$.
It will be convenient to handle the cases $n=7$ and $n=8$ separately.

Let $n=7$.  Then we have
$$
\nu_{7,\theta}(k)=\left\{
\begin{array}{ll}
k+1, & k\leq 3,\\
k+1/4, & k\geq 4,\\
k+\theta, & k=6,
\end{array}
\right.
$$
$$
\sigma_7(k)=\left\{
\begin{array}{ll}
k+1, & k\leq 3,\\
k+1/4, & k\geq 4,
\end{array}
\right.
$$
and 
$$
\mu_7(k)\le\left\{
\begin{array}{ll}
1, & k=2,\\
k+1/4, & k\geq 3.
\end{array}
\right.
$$
Beginning with the case (\ref{final-2}), we see that $2\leq
\max\{\al_2,4\}\leq\al_1$,  and as indicated above we 
may assume that $\al_0+\al_1+\al_2=10$.
If $\al_2=2$ then $\al_0+\al_1=8$, and  it follows that
$E_7(\mal)\leq 3$. This is satisfactory for (\ref{En}).
Alternatively we have $\al_1 \geq 4$ and $\al_2\geq 3$, so that
$E_7(\mal)\leq 3$ in this case also.
The remaining case (\ref{final-4}) is impossible for $n=7$,
since we would have $\al_0=6$, $\al_1 \geq 4$ and $\al_2 \geq 2$, with
$\al_1+\al_2\leq 5$.

We now turn to the case $n=8$, in which setting we have
$$
\nu_{8,\theta}(k)=\left\{
\begin{array}{ll}
k+1, & k\leq 3,\\
k+1/4, & k\geq 4,\\
k+\theta, & k=7,
\end{array}
\right.
$$
and
$$
\sigma_8(k)=\left\{
\begin{array}{ll}
k+1, & k\leq 3,\\
k+1/4, & k\geq 4.
\end{array}
\right.
$$
Moreover we see that $\mu_8(k) =\sigma_8(k)$ if $k \geq 3$, and 
$\mu_8(k)=1$ otherwise.
We begin with the case (\ref{final-2}), for which we have 
$2\leq \max\{\al_2,5\}\leq \al_1$ and $\al_0+\al_1+\al_2=13$.
In view of (\ref{al-triv}) we conclude that $\al_1=5$ or $6$, and at
most one of $\al_0$ or $\al_2$ is less than $4$.
Hence we deduce that $E_8(\mal)\leq 4$ in this case, which is
clearly satisfactory for (\ref{En}).  Finally we must handle the case
in which (\ref{final-4}) holds.  Then $\al_0=7$, and we have the
inequalities $2\leq\max\{\al_2,5\}\leq\al_1$ and $\al_1+\al_2\le 7$.
Thus either $\al_1+\al_2\le 6$, in which case it is clear that 
$E_8(\mal)\leq 4+\theta/3$, or else we must have $\mal=(7,5,2)$.  One easily
deduces that this final possibility is also satisfactory for (\ref{En}).
This completes the treatment of the planes.

\section{Proof of Theorem \ref{6}: $3$-planes}
\lab{3-plane}

Our last task is to consider the case of $3$-planes $M_j\subset
X_3$, for $(d,n) \in \mcal{E}$.  It follows from (\ref{m-range}) that we
may henceforth assume that $n=8$.  It is worth highlighting that this
section would be redundant were we to have a proof of Conjecture
\ref{hope} in the cases $(d,n)=(3,8)$ and $(4,8)$, since then
$\tilde{F}_3(X)$ would be empty.  In the absence of
such a proof there is still work to be done.
We shall essentially employ a blend of the \cc-method and the \cl-method.
For any $a \in X$ let $C_a^3=C_a(\tilde{F}_3(X))$.
As usual the idea will be to estimate the contribution from the
$3$-planes of height $O_\ve(B^{1+\ve})$ according to 
the value of their smallest 
generator, which will necessarily have height $O_\ve(B^{1/4+\ve})$.  
Beginning with the set of smallest
generators $a$ for which the corresponding cone $C_a^3$ has dimension
at most $5$, we obtain the contribution
$$
\ll  \sum_{\colt{a \in
    X(\Q)}{H(a)\ll_{\ve}B^{1/4+\ve}}} 
\#\{x \in C_{a}^3(\Q): H(x)  \leq B\}
$$
to $N_{X_3}(B)$.  Since $\deg C_{a}^3=O(1)$ for each $a \in
X$, an application of part (i) of Lemma~\ref{glass} yields the contribution 
$$
\ll_{\ve}  (B^{1/4+\ve})^{7+\theta}B^{5+1/4+\ve} \ll_{\ve} B^{7+\theta/4+\ve}.
$$
This is clearly satisfactory for Proposition \ref{Xm}.

We now turn to the case in which the cone $C_a^3$ has dimension
at least $6$.  By (\ref{condition}), the dimension is indeed precisely 6.
Lemma \ref{final} now shows that $d=4$, and that the set
$$
G_a=\{T \in \tilde{F}_3(X): a\in T\}
$$
has $\dim G_a=3$.  If $G_a$ were to contain a 3-plane, this would
produce a 4-plane in $X$, contradicting Lemma \ref{fano}.
Hence Lemma \ref{pila'} implies that 
$$
\#\{T \in G_a: H/2<H(T)\leq H\} \ll_\ve H^{3+\ve}.
$$
Moreover Lemmas \ref{triv} and \ref{triv-lin} imply that if $T\in
G_a$ and $H/2<H(T)\leq H$, then $T$
contains $O(B^3+B^4/H)$ rational points of height
at most $B$.  Hence for given $a \in X$ such that 
$C_a^3$ has dimension $6$, there is an overall  contribution  of
\beq\lab{msri}
\ll_\ve \sum_{H \ll_\ve B^{1+\ve}} H^{3+\ve}.\frac{B^{4+\ve}}{H} \ll_\ve 
B^{6+\ve},
\eeq
from the $3$-planes of height  $O_\ve(B^{1+\ve})$ that have smallest
generator $a$.   It remains to sum this over appropriate values of
$a$.

Let $\Psi \subseteq \tilde{F}_3(X)$ be any integral component, so that
in particular $D(\Psi)=X$, and consider the set of lines
$$
\Phi=\{L \in F_1(X): \mbox{$\exists T \in \Psi$ such that $L \subset T$}\}.
$$
Then $\Phi\subseteq F_1(X)$ is an integral subvariety 
that covers $X$, in the sense that $D(\Phi)=X$.
In particular it follows from Lemma \ref{fano} that 
\beq\lab{gower}
\dim \Phi \leq 9,
\eeq
since $(d,n)=(4,8)$.  We then note that
$
D(\Phi_a)=D(\Psi_a)
$
for any $a \in X$.  Thus, for fixed sets $\Psi$ and $\Phi$, 
we are interested in the points
$a \in X$ for which $D(\Phi_a)$ has dimension $6$.  To get a handle
on this set we employ the stratification leading to 
(\ref{s-i}) and (\ref{s-ii}), in conjunction with (\ref{gower}), 
to conclude that there are subvarieties
$$
X=Z_0 \supseteq Z_1 \supseteq Z_2 \supseteq \cdots,
$$
such that $\deg Z_i=O(1)$, with 
\beq\lab{strid}
\dim Z_i \leq 6-i,\quad (i \geq 1),
\eeq
and 
$$
\dim \Phi_a \leq 3+i,\quad (i\ge 0),
$$
for any $a \in Z_i \setminus Z_{i+1}$, 
for $i\geq 0$.

It now follows that if $\dim D(\Phi_a)=6$ then $\dim\Phi_a=5$, and
hence that $a\in Z_i$ for some $i\ge 2$.  Moreover we will have $\dim
Z_i\le 4$.  Thus the overall contribution to $N_{X_3}(B)$
corresponding to points with $H(a)\le A$ will be
$$
\ll_\ve A^{4+1/4+\ve}B^{6+\ve} 
$$
in view of part (i) of Lemma \ref{glass} and the estimate (\ref{msri}).  On
choosing $A=B^{4/17}$ we see that this gives a satisfactory treatment
of $3$-planes whose smallest generator has height at most $B^{4/17}$.

It now remains to deal with the case in which the smallest generator
$a$ of $T$ has height larger than $B^{4/17}$.
Now it may obviously be assumed that none of the
$3$-planes in which we are interested lie in $Z_1$.  Indeed Lemma
\ref{triv} implies that there are at most $O(B^6)$ points of height at
most $B$ contained in $Z_1$, which is satisfactory for
Proposition~\ref{Xm}.  
We have already seen in \S \ref{campaign} that any
$3$-plane $T$ defined over $\Q$ contains linearly 
independent points $a_0,a_1,a_2,a_3 \in \bfP^8(\Q)$ such that
$$
H(a_0) \le \cdots \le H(a_3), \quad H(T) \ll H(a_0)H(a_1)H(a_2)H(a_3) \ll H(T).
$$
Since $H(a_0)\geq B^{4/17}$ and $H(T)\ll_{\ve} B^{1+\ve}$ it follows that
$H(a_3)\ll_{\ve} B^{5/17+\ve}$.  In particular $H(a_i)\ll B^{1/3}$ for
$i=0,\ldots,3$.   We claim that there is a point $b\in T\setminus Z_1$
for which $H(b)\ll B^{1/3}$.  Indeed there are $\gg L^4$ points
\[
b=\lambda_0a_0+\lambda_1a_1+\lambda_2a_2+\lambda_3a_3\in T,\quad
1\le\lambda_i\le L,
\]
of which $O(L^3)$ can lie on the proper subvariety $T\cap Z_1$, by
Lemma \ref{triv}.  Thus if $L$ is a sufficiently large constant we
may produce at least one point $b\in T\setminus Z_1$, with $H(b)\ll
B^{1/3}$.

For each 3-plane $T$ under consideration we now choose an appropriate
point $b$.  We then proceed to count points on $T$ according to the
corresponding values of $b$, noting that if $x\in T$ then $x\in C_b(\Psi)=D(\Psi_b)$.
Since $b\not\in Z_1$ we will have 
\[\dim C_b(\Psi)\le \dim C_b^1 = \dim D(\Phi_b)= 1+\dim \Phi_b\le 4,\]
by (\ref{strid}).  According to part (i) of Lemma \ref{glass} the total
contribution to $N_{X_3}(B)$ arising in this way is then
$$
\leq  \sum_{\colt{b \in
    X(\Q)}{H(b)\leq B^{1/3}}} \#\{x \in C_{b}(\Psi)(\Q): H(x)  \leq B\}
\ll_\ve (B^{1/3})^{7+\theta}B^{4+1/4+\ve},
$$
which is satisfactory for Proposition \ref{Xm}.
This completes the proof of 
Proposition~\ref{Xm} in the case $m=3$.

\newpage
\appendix

\section*{
\begin{center}
Appendix
\end{center}}

\begin{center}
{J.M. Starr}\\
\small{{\em Department of Mathematics, MIT, Cambridge MA 02139}}\\
\small{jstarr@math.mit.edu} 
\end{center}

\begin{abstract}
A smooth, nondegenerate
hypersurface in projective space contains no linear subvarieties of
greater than half its dimension.  It can contain
linear subvarieties of half its dimension.  This note proves that a
smooth hypersurface of degree $d\geq 3$
contains at most finitely many such subvarieties.
\end{abstract}

\medskip\noindent 
Let $k$ be a field.  Let $X\subset \PP^n_k$ be a hypersurface of
degree $d > 1$.  For each integer $m > 0$, denote by $F_m(X)$ the
\emph{Fano scheme of $m$-planes in} $X$, cf. \cite{AKFano}.  There
are a number of natural questions about $F_m(X)$: is this scheme
non-empty, is this scheme connected, is this scheme irreducible, is
this scheme reduced, is this scheme smooth, what is the dimension of
this scheme, 
what is the degree of this
scheme, etc.?  For each of these questions, the answer is uniform for
a \emph{generic hypersurface}.  More precisely, 
there is a non-empty open subset $U$ of the parameter space 
of hypersurfaces, such that for
every point in $U$ the answer to the question is the same.  For a generic
hypersurface, the answer is often easy to find:  the total space of the
relative Fano scheme of the universal hypersurface is itself a
projective bundle over the Grassmannian $\mathbb{G}(m,n)$, so if 
the question for a generic hypersurface can be reformulated as a
question about the total space of the relative Fano scheme, it is easy
to answer the question.
However, much less
is known if $X$ is assumed to be smooth, but \emph{not} generic.

There are a few easy results, such as the following.

\begin{prop} \label{prop-1}
Let $X \subset \PP^n$ be a hypersurface of degree $d > 1$ and let $m$
be an integer such that $2m \geq n$.
Every $m$-plane $\Lambda \subset X$ intersects the singular locus of $X$.
In particular, if $X$ is smooth then $F_m(X)$ is empty.
\end{prop}

\begin{proof}
Choose a system of homogeneous coordinates $x_0,\dots,x_n$ on $\PP^n$
such that $\Lambda$ is given by $x_{m+1}=\cdots =x_n=0$.  
Let $F$ be a defining equation for $X$.  Because $\Lambda \subset X$,
$F(x_0,\dots,x_m,0,\dots,0) = 0$.  Also,
$$
\frac{\partial F}{\partial x_i}(x_0,\dots,x_m,0,\dots,0) = 0,
$$
for $i=0,\dots,m$.
Because $d>1$, for $i=1,\dots,n-m$, the homogeneous polynomial on $\Lambda$,
$$
\frac{\partial F}{\partial x_{m+i}}(x_0,\dots,x_m,0,\dots,0),
$$
is non-constant.  Since $n-m \leq m$, these $n-m$ non-constant
homogeneous polynomials 
have a common zero in $\Lambda$. 
By the Jacobian criterion, this is a singular point of $X$.
\end{proof}

\begin{rmk}  This also follows easily from the Lefschetz hyperplane
  theorem.
\end{rmk}

What happens if $n=2m+1$?  If $m>1$ or if $m=1$ and $d>3$, then a
\emph{generic} hypersurface $X\subset \PP^n$ contains no $m$-plane.
However there do exist smooth hypersurfaces 
containing an $m$-plane.  For instance, 
if $\text{char}(k)$ does not divide
$d$ then the Fermat hypersurface $x_0^d + \dots + x_n^d=0$ 
is smooth and contains many $m$-planes, e.g.,
$x_0+x_1=x_2+x_3=\cdots= x_{n-1}+x_n=0$ when $d$ is odd.  
However, if $d \geq 3$,
a smooth
hypersurface cannot contain a positive-dimensional family of
$m$-planes.  This was proved independently by Olivier Debarre, using a
different argument.  

The setup is as follows.  
Let $n=2m+1$.  Let $X \subset \PP^n$ be a hypersurface of
degree $d$.  Let $\Lambda_1,\Lambda_2 \subset X$ be $m$-planes.
Denote by $Z$ the intersection $\Lambda_1 \cap \Lambda_2$.  This is
either
empty or else an
$r$-plane for some integer $r$.  If $Z$ is empty, define $r$ to be
$-1$.

Denote by $X_\text{sm} \subset X$ the smooth locus of $X$, i.e., the
maximal open subscheme that is smooth.  Denote $\Lambda_{i,\text{sm}}
= \Lambda_i \cap X_\text{sm}$ for $i=1,2$.
There are Chow classes,
$$
[\Lambda_{i,\text{sm}}] \in A_m(X_\text{sm}), \
i=1,2.
$$  
Because
$X_\text{sm}$ is smooth, the intersection product
$[\Lambda_{1,\text{sm}}] \cdot 
[\Lambda_{2,\text{sm}}] \in A_0(X_\text{sm})$ is defined.  

\begin{lemma} \label{lem-2}
If $Z$ is contained in $X_\text{sm}$, then the degree of
$[\Lambda_{1,\text{sm}}] \cdot
[\Lambda_{2,\text{sm}}]$ is $( 1 - (1-d)^{r+1})/d$. 
\end{lemma}

\begin{proof}
If $r=-1$, i.e., if $Z$ is empty, this is obvious.  Therefore suppose
that $Z$ is an $r$-plane for some $r \geq 0$.  By the excess
intersection formula, the class $[\Lambda_{1,\text{sm}}]\cdot
[\Lambda_{2,\text{sm}}]$ is the
pushforward from $Z$ of the refined intersection product,
$(\Lambda_{1,\text{sm}}\cdot \Lambda_{2,\text{sm}})^Z$.  And by
\cite[Prop. 9.1.1]{F}, the refined intersection product is,
$$
(\Lambda_{1,\text{sm}}\cdot \Lambda_{2,\text{sm}})^Z = \lt\{
c(N_{\Lambda_{1,\text{sm}}/X_{\text{sm}}})/c(N_{Z/\Lambda_{2,\text{sm}}})\cap
[Z] \rt\}_r.
$$
Denote $H = c_1(\OO_Z(1))$.  The normal bundle of $Z$ in
$\Lambda_{1,\text{sm}}$ is $\OO_Z(1)^{m-r}$.  The restriction to $Z$
of the normal bundle of $\Lambda_{2,\text{sm}}$ in $\PP^n$ is
$\OO_Z(1)^{n-m}= \OO_Z(1)^{m+1}$.  
And the restriction to $Z$ of the normal bundle of
$X$ in $\PP^n$ is $\OO_Z(d)$.  Therefore,
$$
c(N_{\Lambda_{1,\text{sm}}/X_{\text{sm}}})/c(N_{Z/\Lambda_{2,\text{sm}}})
= \frac{(1+H)^{m+1}}{(1+H)^{m-r}(1+dH)} = \frac{(1+H)^{r+1}}{1+dH}.
$$
Expanding this out gives,
$$
\left(\sum_{i=0}^{r+1} \binom{r+1}{i} H^i\right) 
\left(\sum_{j=0}^{\infty} (-1)^j 
  d^j H^j\right). 
$$
In particular, the coefficient of $H^r$ is,
\begin{align*} 
\sum_{i=0}^r \binom{r+1}{i} (-1)^{r-i} d^{r-i} &= 
\frac{-1}{d} \sum_{i=0}^r \binom{r+1}{i} (-1)^{r+1-i} d^{r+1-i}\\ 
&= ( 1 - (1-d)^{r+1})/d.
\end{align*}
\end{proof}

\begin{prop} \label{prop-3}
If $d\geq 3$, if $\Lambda_1$ and $\Lambda_2$ are distinct, and if at
least one of $\Lambda_1$, $\Lambda_2$ is contained in $X_{\text{sm}}$,
then $[\Lambda_{1,\text{sm}}]$ is not numerically equivalent to
$[\Lambda_{2,\text{sm}}]$.  
\end{prop}

\begin{proof}
Let $\Lambda_1$ be contained in $X_{\text{sm}}$.  By
Lemma \ref{lem-2},
$$
\text{deg}([\Lambda_1]\cdot [\Lambda_1]) = ( 1 - (1-d)^{m+1})/d. 
$$
Also, $Z = \Lambda_1\cap \Lambda_2$ is contained in $X_\text{sm}$.  So
by Lemma~\ref{lem-2},
$$
\text{deg}([\Lambda_1]\cdot[\Lambda_{2,\text{sm}}]) =
( 1 - (1-d)^{r+1})/d. 
$$
Because $d-1 \geq 2$ and $r < m$, we have  
$(d-1)^{r+1} < (d-1)^{m+1}$. Therefore $( 1 - (1-d)^{r+1})/d\neq  ( 1 - (1-d)^{m+1})/d$, 
and so $[\Lambda_1]$ is not numerically equivalent to 
$[\Lambda_{2,\text{sm}}]$.
\end{proof}

\begin{coro}[Debarre] 
There are only finitely many $m$-planes contained in $X_{\text{sm}}$.
\end{coro}

\begin{proof}
By Proposition \ref{prop-3}, distinct $m$-planes contained in
$X_{\text{sm}}$ are not algebraically equivalent.  Therefore every
irreducible component of $F_m(X)$ that contains a point parametrising
an $m$-plane in $X_{\text{sm}}$ is just a point.  Because $F_m(X)$ is
quasi-compact, the number of these irreducible components is finite.  
\end{proof}

\begin{rmk} \label{rmk-5}
Debarre's proof shows more than the statement of the corollary:
for any $m$-plane $\Lambda$ contained in
$X_{\text{sm}}$, $h^0(\Lambda,N_{\Lambda/X}) = 0$.  It follows that
each such point is a connected component of $F_m(X)$, and that
$F_m(X)$ is reduced at this point.  
\end{rmk}

A natural question is, what is the maximal number of $m$-planes
contained in a smooth hypersurface of degree $d$ in $\PP^{2m+1}$?  There is
a naive upper bound that grows as $d^{(m+1)^2}$, but this is
too large.  The Fermat hypersurface contains $C_m d^{m+1}$
distinct $m$-planes, where $C_m = (2m+1)(2m-1)\cdot \dots \cdot 3
\cdot 1$.  Joe Harris points out that for $m=1$, the degree of the
\emph{flecnodal curve} gives an upper bound of $11d^2-24d$.  

More
generally, define the \emph{flecnodal locus} $P(X) \subset X$ to be the set of
points $p\in X$ such that there is a line $L \subset \PP^{2m+1}$ that
has contact of order $3m+1$ with $X$ at $p$.  This is the pushforward
in $X$ of a subscheme in $\PP(T_X)$ that is the zero locus of a
section of a locally free sheaf.  If $X$ is generic, this section is a
regular section.  Then a Chern class computation gives that the degree
of $P(X)$ is a polynomial $p_m(d)$ of degree $m+1$ in $d$ whose leading term
is,
$$
\lt( \frac{(3m+1)!}{2} - 1 \rt) d^{m+1}.
$$
Therefore, for arbitrary $X$, $\text{deg}(P_m(X)) \leq p_m(d)$, where
$P_m(X)$ is the $m$-cycle of all $m$-dimensional irreducible
components of $P(X)$ (weighted by multiplicity).

Of course every
$m$-plane $\Lambda$ is contained in $P(X)$.  It is not clear that
every $m$-plane is contained in $P_m(X)$, i.e., that $\Lambda$ is an 
irreducible component of $P(X)$.  And, indeed, this
fails if $d \leq 3m$.  For $d \geq 3m$, it may be true that every
$m$-plane is an irreducible component of $P(X)$.

In the special case that $m=1$, 
$P(X)$ is a curve, the \emph{flecnodal curve}, for all $d\geq 3$.
Therefore, the number of lines in a smooth surface of
degree $d\geq 3$ in $\PP^3$ is at most
$11d^2-24d$.  Note this gives the correct answer for $d=3$.  For $d=4$
this gives the wrong answer,
Segre proved the maximal number of lines on a quartic surface is $64$,
cf. \cite{SegreLines}.  In fact, by a more involved analysis, Segre
proved that the number of lines on a smooth surface of degree $d\geq 3$
is at most $11d^2 -28d+12$.
The Fermat surface
contains $3d^2$ lines.  The true maximum is probably strictly 
between $3d^2$ and $11d^2-28d+12$.

\end{document}